\documentclass{article}
\usepackage{graphicx}
\usepackage{verbatim}
\usepackage{amsmath,amsfonts,amssymb,amsthm}
\usepackage{latexsym}
\usepackage{multirow}
\usepackage{fancyhdr}
\usepackage{amsmath,amscd}
\usepackage{enumitem}
\usepackage{float}
\theoremstyle{plain}
\newtheorem{theorem}{Theorem}[section]
\newtheorem{theorem1}{Theorem}
\newtheorem{proposition}[theorem]{Proposition}
\newtheorem{lemma}[theorem]{Lemma}
\newtheorem{corollary}[theorem]{Corollary}

\theoremstyle{definition}
\newtheorem{definition}[theorem]{Definition}
\newtheorem{example}[theorem]{Example}
\theoremstyle{remark}
\newtheorem{remark}[theorem]{Remark}
\theoremstyle{remark}

\theoremstyle{remark}

\theoremstyle{remark}

\providecommand{\keywords}[1]
{
	\small	
	\textbf{Keywords:} #1
}

\providecommand{\subjclass}[1]
{
	\small	
	\textbf{MSC-Class:} #1
}
  
\title{Conjugacy classes of groups of prime order in $\mathrm{PGL}_{k+1}(\mathbb{C})$}
\author{Andrea Marinatto\\
 e-mail: ettabon@alice.it}

\begin{document}

\maketitle

\begin{abstract}
Let $\mathbb{C}$ be the field of complex numbers. Let $k$ be natural number with $k \geq 2$ and let $p$ be a rational prime. In this paper we count the number of conjugacy classes of admissible cyclic subgroups  of $\mathrm{PGL}_{k+1}(\mathbb{C})$ of order $p$, where with admissible we intend those finite subgroups that can be contained in the automorphism group of a set of points  in  $\mathbb{P}^k(\mathbb{C})$ in general position and of cardinality $n\geq k+3$. We also describe a kind of association between the conjugacy classes of these groups and show  a beautiful relation connecting this type of association and the association between point sets.
\end{abstract}



\section{Introduction}
	Let $\mathbb{C}$ be the field of complex numbers. The aim of this  paper is to count the number of conjugacy classes of certain admissible subgroups of   $\mathrm{PGL}_{k+1}(\mathbb{C})$ mainly for values of  $k \in \left\lbrace {2,3,4,5}\right\rbrace$. With admissible subgroups  we intend those finite subgroups of $\mathrm{PGL}_{k+1}(\mathbb{C})$ that can be contained in the automorphism group of a set of points  in  $\mathbb{P}^k(\mathbb{C})$ in general position and of cardinality $n\geq k+3$. The concept of admissible 	subgroups originated in in~\cite{Mari2}, where we needed to describe and classify the admissible subgroups of prime order to bound from the above the dimension of the isomorphism
	classes of point sets in $\mathbb{P}^k(\mathbb{C})$ whose automorphism is not trivial. We carried out that classification by determining the pattern that the eigenvalues of a matrix of finite order 	$p$ have to follow for this matrix to represent an admissible automorphim in $\mathrm{PGL}_{k+1}(\mathbb{C})$. 	In that paper we also addressed the question whether for an abstract finite group there 	is more than one admissible group (up to conjugation) which realizes it as a subgroup 	of $\mathrm{PGL}_{k+1}(\mathbb{C})$. This analysis has led us to the computation of the number of conjugacy classes of admissible group we are carrying out here. The results we show in this paper cover two areas of interest. The first one is naturally the study of the finite subgroups
	of $\mathrm{PGL}_{k+1}(\mathbb{C})$. Apart for the computation of the conjugacy classes of admissible cyclic groups of order $p$, our analysis gives information on more complex groups that	can be obtained as semidirect products of one of these cyclic groups and a group generated by the automorphism represented by a permutation matrix, thus partially extending the analysis given in~\cite{Dolgar} for the finite subgroups of  $\mathrm{PGL}_{3}(\mathbb{C})$: this follows from our description of the cases where the	normalizer of an admissible cyclic group of order $p$ contains as subgroup the group generated by
	the automorphism represented by a permutation matrix. The second application of our results is in the study of the point sets and their automorphisms. Point sets play an important role in several branches of algebraic geometry. For instance the isomorphism classes of $n$-point set (for $n\geq 4$) on  are related to the isomorphism classes of
	elliptic and hyperelleptic curves and binary forms. Except for the case of $4$ points on $\mathbb{P}^1(\mathbb{C})$, a generic point set has trivial automorphism group so the cases where a point set has extra automorphisms are of special interest. For example, when dealing with the rationality problem of the equivalence between fields of moduli and field of definition for point sets on $\mathbb{P}^1(\mathbb{C})$ in~\cite{Mari1}, the cohomological approach has to be modified to treat the non-trivial case so requiring additional analysis. Our contribution in this paper as well in~\cite{Mari2} is to show which cyclic group of prime order can be contained in the	automorphims group of a point set and also describe a kind of canonical form for these point sets.\\
	With regard to the contents of this paper, the computation of the number of the 	conjugacy classes of admissible subgroups of order $p$ in $\mathrm{PGL}_{k+1}(\mathbb{C})$ requires to split the analysis into two subcases, $p\geq  k + 1$ and $p < k + 1$, and our treatment of the subject proceeds through different steps.		
In section~\ref{section0}  we recall some basic notions regarding the sets of points in general positions in $\mathbb{P}^{k}(\mathbb{C})$ and the admissible cyclic subgroups of $\mathrm{PGL}_{k+1}(\mathbb{C})$. Here we recall in particular that a finite cyclic  subgroup of $\mathrm{PGL}_{k+1}(\mathbb{C})$ can be put in diagonal form and that when in this shape and has order $m$ any of its generators can be represented  by a diagonal matrix  in  $\mathrm{GL}_{k+1}(\mathbb{C})$ of finite order $m$,  whose entries are all powers of a  primitive $m$-th root of unity. Then, we restrict our attention to the cyclic groups of prime order $p\geq k+1$ and recall that they are admissible if and only if the entries of this diagonal matrix are all distinct. In this section we also remind some definitions and key facts regarding  the sets of points in general positions.  Then, in Section~\ref{section1} we achieve the computation of the number of conjugacy classes of these admissible subgroups. The computations proceeds through different steps. We first note that each of these groups can be uniquely identified by a matrix of finite order $p$ and of assigned shape. In the paper this matrix is referred to  as the identifying matrix. The first step in our analysis is then to compute the total number of these matrices. This is equivalent to compute the  number of certain $(k+1)$-tuples of integers modulo the prime $p$.  This $(k+1)$-tuples corresponding to  the identifying matrices are described as identifying vectors. Then we consider the action of the symmetric group $S_{k+1}$ on the set of  the subgroups of $\mathrm{PGL}_{k+1}(\mathbb{C})$  of order $p$ and in diagonal form. Recall that the symmetric group $S_{k+1}$ can be identified as the group of permutation matrices in $\mathrm{PGL}_{k+1}(\mathbb{C})$. If  $H$ is the stabilizer  of an admissible group $G$ for the action of $S_{k+1}$ then the number of its diagonal conjugates equals the index $\left[S_{k+1}:H\right]$ of $H$ in $S_{k+1}$. This can be read in terms of identifying matrices.  This leads to reverse the problem,  that is fixing a subgroup $H$ in $S_{k+1}$ to see how many identifying matrices (if there are any) generate a subgroup whose stabilizer is $H$ . For the cyclic subgroups of $S_{k+1}$ in particular,  this counting is made by considering the entries of an identifying  matrix as unknowns and  consider the permutation on these entries induced by a generator of  the group. This results in a system of congruences modulo $p$. If this system has admissible solution  (i.e. tuples whose entries are all distinct),  these solutions are the entries of the identifying matrix with the requested condition. For our purposes we should consider all the subgroups belonging to $S_{k+1}$, but our examination is highly simplified by the fact that the number of solutions is constant along the conjugacy classes. Moreover the complete description of the lattice of  subgroups of $S_{k+1}$ for $k \in \left\lbrace {2,3,4,5}\right\rbrace$ allows to rapidly reject subgroups containing permutations whose systems of congruences are not admissible. The outcome  of our examinations is that the  existence of solutions for the subgroups generated by cyclic permutations of  $S_{k+1}$ of length $k$ or $k+1$  depends on the parity of the prime $p$ modulo certain integers linked to these permutations.  On the other hand, subgroups generated by cycles of smaller length or more  structured subgroups, such for example the dihedral ones, do not  lead  to  admissible  systems of congruences except in a very few cases.  In our analysis, we also check for any possible overlap of solutions between systems pertaining to distinct subgroups.  We have found  only few overlaps and almost all  of them  are between systems relative to subgroups of a cyclic group. In this case the solutions in common have to be taken out from the counting of those relative to the smaller subgroups. \\
Then we pass to the computation of the number of conjugacy classes. We fix a prime $p$  and a non-trivial subgroup $H$ of $S_{k+1}$  whose system of congruences modulo $p$ is  admissible and count the solutions. Then we multiply this number by the length of the conjugacy class of $H$. This product divided by the index of $H$ in $S_{k+1}$ is the number of admissible conjugacy classes of subgroups of order $p$ whose stabilizer belongs to the class of $H$. Then we take care of  the conjugacy classes of admissible group of order $p$ whose stabilizer is trivial. This is done by subtracting from the  total number of identifying matrices the number of solutions relative to the non-trivial subgroups of $S_{k+1}$ and dividing the number so obtained by the factorial of $k+1$. We have applied this method for $k \in \left\lbrace {2,3,4,5}\right\rbrace$  and our results are listed here below in Theorems from $1$ to $5$. The number $\Phi$ appearing in the tables below stands for the total number of identifying matrices. It obviously depends on $k$ and $p$ and is explicit expression is given by  $\Phi=(p-2)\cdots (p-k+1)$. Our  first results regard the case $k=2$ and  are listed in Theorem  $1$ below. This is the easiest case we have treated  due to the small lattice of subgroups of $S_3$ and the fact  there is no overlap of solutions. 
\setcounter{theorem1} {0} 
\begin{theorem1}(see Theorem~\ref{1theorem} )
		Let  $p>3$ be a prime. Then the number of conjugacy classes of admissible cyclic groups of order $p$ in $\mathrm{PGL}_{3}(\mathbb{C})$ is as follows depending on the parity modulo $3$ of $p$:
		\begin{table}[H]
			\centering
					\renewcommand{\arraystretch}{1.5}
			\begin{tabular}{|c|c|}
								\hline 
					Prime Parity &  Number of Conjugacy Classes \\ 
				\hline 
				  $p\not\equiv 1 \ (mod \ 3)$ & $1+ \frac{\Phi-3}{6}$ \\
				\hline 
				   $p\equiv 1 \ (mod \ 3)$ &  $2+  \frac{\Phi-5}{6}$ \\ 
				\hline 
			\end{tabular} 
		\end{table}
\end{theorem1}	
The next case we have treated is $k=3$, and the relative results are listed in Theorem $2$. This value of $k$ marks the first occurrence of overlaps between solutions. More specifically we have found some common solutions between the systems relative to the $4$-cycles of $S_4$ and the permutations which are the products of two distinct two  cycles.
\begin{theorem1}{(See Theorem~\ref{acase3})}
			Let $p>3$ be a prime. Then the number of conjugacy classes of admissible cyclic groups of order $p$ in $\mathrm{PGL}_{4}(\mathbb{C})$ is as follows depending on the parity of $p$ modulo $3$ and $4$: 
		\begin{table}[H]
			\centering
			\renewcommand{\arraystretch}{2}
			\begin{tabular}{|c|c|}
				\hline 
					Prime Parity &  Number of Conjugacy Classes \\ 
				\hline 
					$p\not\equiv 1 \ (mod \ 4)\wedge p\not\equiv 1 \ (mod \ 3)$ & $\frac{\alpha}{12} +\frac{\Phi-\alpha}{24}$  \\ 
				\hline 
				
					$p\not\equiv 1 \ (mod \ 4)\wedge p\equiv 1 \ (mod \ 3)$  & $1+\frac{\alpha}{12} +\frac{\Phi-\alpha-8}{24} $\\ 
				\hline  
					$p\equiv 1 \ (mod \ 4)\wedge p\not\equiv 1 \ (mod \ 3)$ & $1+\frac{\alpha-6}{12} +\frac{\Phi- \alpha}{24} $  \\ 
				\hline 
					$p\equiv 1 \ (mod \ 4)\wedge p\equiv 1 \ (mod \ 3)$ & $2+\frac{\alpha-6}{12} +\frac{\Phi- \alpha-8}{24}  $ \\ 
				\hline
				\multicolumn{2}{l} {where $\alpha=3\cdot (p-3)$.}
			\end{tabular} 
		\end{table}
\end{theorem1} 
The results for the last two cases treated, $k=4$ and $k=5$, are displayed in Theorems $3$	 and $4$. The computations here  required more effort due to the complexity of the lattices of subgroups of $S_5$ and $S_6$.  Here again there is the presence of overlaps.

\begin{theorem1}{(See Theorem~\ref{case4})}
	Let $p>5$ be a prime. Then the number of conjugacy classes of admissible cyclic groups of order $p$ in $\mathrm{PGL}_{5}(\mathbb{C})$ is as follow depending on the parity of $p$ modulo $4$ and $5$:
	\begin{table}[H]
		\centering
		\renewcommand{\arraystretch}{2}
		\begin{tabular}{|c|c|}
			\hline 
				Prime Parity &  Number of Conjugacy Classes \\ 
			\hline	
			 $p\not\equiv 1 \ (mod \ 5)\wedge p\not\equiv 1 \ (mod \ 4)$ & $\frac{\beta}{60} +\frac{\Phi-\beta}{120}$  \\ 
			\hline 	
			$p\not\equiv 1 \ (mod \ 5)\wedge p\equiv 1 \ (mod \ 4)$  & $1+\frac{\beta-30}{60} +\frac{\Phi-\beta}{120}$\\ 
			\hline  
			 $p\equiv 1 \ (mod \ 5)\wedge p\not\equiv 1 \ (mod \ 4)$ & $1+\frac{\beta}{60} +\frac{\Phi-\beta-24}{120}$  \\ 
			\hline 
			 $p\equiv 1 \ (mod \ 5)\wedge p\equiv 1 \ (mod \ 4)$ & $2+\frac{\beta-30}{60} +\frac{\Phi-\beta-24}{120}$ \\ 
			\hline 
			\multicolumn{2}{l} {where $\beta=15\cdot (p-3)$.}
		\end{tabular} 
	\end{table}
\end{theorem1}

\begin{theorem1}(see Theorem~\ref{k5})
		Let $p\geq 7$  be a prime. Then, the number of conjugacy classes of admissible cyclic subgroups of order $p$ in $\mathrm{PGL}_{6}(\mathbb{C})$ is as follows depending on the parity of $p$ modulo  $5$ and $3$:
		\begin{table}[H]
						\begin{center}
				\renewcommand{\arraystretch}{2}			
				\begin{tabular}{|c|c|}\hline
					Prime Parity & Number of Conjugacy Classes \\ \hline
					 $p\not\equiv 1 \ (mod \ 5)\wedge p\not\equiv 1 \ (mod \ 3)$ & $\frac{\alpha}{360}+ \frac{\Phi-\alpha}{720}$\\ \hline
					 $p\equiv 1 \ (mod \ 5)\wedge p\not\equiv 1 \ (mod \ 3)$ & $1+\frac{\alpha}{360} +\frac{\Phi-144-\alpha}{720}$ \\ \hline
					 $p\not\equiv 1 \ (mod \ 5)\wedge p\equiv 1 \ (mod \ 3)$ & $1+\frac{\beta}{360}+\frac{\delta}{240}+\frac{\Phi-120-\beta-\delta}{720}$ \\
					\hline
					 $p\equiv 1 \ (mod \ 5)\wedge p\equiv 1 \ (mod \ 3)$ & $2+\frac{\beta}{360}+\frac{\delta}{240} +\frac{\Phi-120-144- \beta-\delta }{720}$	\\ \hline 
					 	\multicolumn{2}{l} {where $\alpha=(p-3)\cdot(p-5)\cdot 15$, $\beta=(p-3)\cdot(p-5)\cdot 15-120$ and $\delta=2\cdot (p-4)\cdot 20-120$ }
				\end{tabular}
			\end{center}
			
		\end{table} 
\end{theorem1}
The results for the special case $p=k+1$ are not listed in Theorem $1$ and in Theorem $3$. For these values of $p$ there is always a single conjugacy class.  \\
In Section~\ref{othergroups} we show that the results given in Section~\ref{section1} can be used to compute the conjugacy classes of certain other admissible groups that	are semidirect products of  one of the  admissible cyclic groups we have described so far and a group generated by some permutation matrix.
In Section~\ref{power} we have treated the case  $p<k+1$. Here, our main result is the following theorem
\begin{theorem1}(see Theorem~\ref{aconjugacy})\label{dualitytheorem}
	Let $p\geq 5$ be a prime with $p<k+1$. Suppose further that $k+1\equiv a \ (mod \ p)$ with $2\leq a \leq p-2$.  Then, the number of conjugacy classes of  admissible subgroups of order $p$ in $\mathrm{PGL}_{k+1}(\mathbb{C})$ is the same as the number of conjugacy classes of  admissible subgroups of order $p$  in $\mathrm{PGL}_{a}(\mathbb{C})$ and in $\mathrm{PGL}_{p-a}(\mathbb{C})$.
\end{theorem1}
This theorem  has two application in the paper. The first is to compute  the  number of conjugacy classes of admissible subgroups of order $p$ when $p<k+1$ using the results given in theorems from $1$ to $5$ above for $p\geq k+1$. This application, together with the  results for the cases $k+1\equiv 0, 1 \ (mod \ p)$,  permits in particular to obtain the number of  conjugacy classes of the primes in the set $\lbrace{2,3,5,7,11,13}\rbrace$ for any value of $k$.   The  second use of Theorem~\ref{dualitytheorem} exploits the kind of duality between $\mathrm{PGL}_{a}(\mathbb{C})$ and in $\mathrm{PGL}_{p-a}(\mathbb{C})$. In this way we could partially extend  the results  given in the previous theorems for $p\geq k+1$ to  projective general linear groups of degree higher than $6$.  As an example we have computed the number of conjugacy classes of subgroups of order $19$ for values of $k+1$  ranging from $13$ and $19$. 
In this section we also show a beautiful relation connecting  the association between cyclic groups of order $p$ and the association between point sets of cardinality $p$ consisting of one non-trivial orbit under the action of these groups. This is the content of the following theorem
\begin{theorem1}{(see Theorem~\ref{correspondance})}
	Let $p\geq 5$ be a prime. Suppose further that $a$ is an integer with  $2\leq a \leq p-2$.  Let $H$ be an  admissible subgroup of order $p$ in $\mathrm{PGL}_{a}(\mathbb{C})$ in diagonal form and let $K$ be one of its associated subgroups in $\mathrm{PGL}_{p-a}(\mathbb{C})$ described in Theorem~\ref{aconjugacy}.
	Let $T$ be the point set in $\mathbb{P}^{a-1}\left(\mathbb{C}\right)$ consisting of the orbit of the identity point under the action of the group $H$. Similarly, let $V$ be the point set in $\mathbb{P}^{p-a-1}\left(\mathbb{C}\right)$ consisting of the orbit of the identity point under the action of the group $K$. Then the point sets $T$ and $V$ are associated.
\end{theorem1}
We also show how to extend the association described in the above theorem to point sets consisting of more than one non-trivial orbit.
\section{Background Notions}\label{section0}

In this  section  we introduce the basic notions regarding point sets and  admissible cyclic subgroups of $\mathrm{PGL}_{k+1}(\mathbb{C})$. We start with those the regarding point sets

\subsection{Point Sets}
In the remainder of this paper $\mathbb{C}$ is the the field of complex numbers. Moreover,  $n$ and $k$ are natural numbers satisfying $k \geq 2$ and $n\geqslant k+3$ and $\mathbb{P}^k_n\left(\mathbb{C}\right)$ be the subset of $n$-th symmetric product $(\mathbb{P}^k(\mathbb{C}))^{(n)}$ of $\mathbb{P}^k(\mathbb{C})$ with itself consisting of $n$-point sets in general position. In this paper the elements of $\mathbb{P}^k_n\left(\mathbb{C}\right)$ will be referred to simply as $n$-point sets. \\

The action of the group $\mathrm{PGL}_{k+1}\left(\mathbb{C}\right)$ on $\mathbb{P}^k\left(\mathbb{C}\right)$ induces an action on $\mathbb{P}^{k}_n\left(\mathbb{C}\right)$ given for each $T=\lbrace P_1, \dots, P_n\rbrace \in \mathbb{P}^{k}_n(\mathbb{C})$ and each $f\in \mathrm{PGL}_{k+1}\left(\mathbb{C}\right)$ by 
\begin{center}
	$T^f= \left\{ P_1^f, \ldots, P_n^f\right\}$.
\end{center}

The orbit of a point set $T$ under the action of $\mathrm{PGL}_{k+1}\left(\mathbb{C}\right)$ will be denoted with $\left[T\right]$: in other words $\left[T\right]$ is the image of $T$ in the quotient space

\begin{center}
	$M_n^{k}(\mathbb{C})\stackrel{def}{=} \mathbb{P}^{k}_n\left(\mathbb{C}\right)/\mathrm{PGL}_{k+1}\left(\mathbb{C}\right)$.
\end{center} 
The quotient set $M_n^k(\mathbb{C})$ admits a natural structure of algebraic variety such that the canonical projection from $\mathbb{P}^k_n\left(\mathbb{C}\right)$ is a morphism (see for reference Chapters $5$, $6$ and $7$ in~\cite{Mukai}). Note also  that the dimension of $M_n^k(\mathbb{C})$ is $n\cdot k-k\cdot(k+2)$.\\

For each $T\in \mathbb{P}^k_n$, the group of automorphisms of $T$ is the group
\begin{center}
	$\mathrm{Aut}\left(T\right)=\left\{f\in \mathrm{PGL}_{k+1}\left(\mathbb{C}\right)|T^f=T\right\}$.
\end{center}.
\begin{proposition}\label{groups}
	Let let $n$ and $k$ be positive integers satisfying $n\geqslant k+3$ and let $T\in \mathbb{P}_n^k\left(\mathbb{C}\right)$ be a $n$-point set, then
	
	\begin{description}
		\item [1)] $\mathrm{Aut}(T)$ is a finite group.
		\item [2)] For any $f\in \mathrm{PGL}_{k+1}\left(\mathbb{C}\right)$ the map $\phi$
		
		\begin{center}
			$\phi:\mathrm{Aut}\left(T\right)\rightarrow \mathrm{Aut}\left(T^f\right)$\\
			
			$h\mapsto  f\circ h\circ f^{-1}$ 
		\end{center}
		is an isomorphism of groups.  
	\end{description}
\end{proposition}

\begin{proof}
	\begin{description}
		\item [(a)] With the necessary modifications the proof is the same as that of (a) in Prop 2.3 in~\cite{Mari1}.
		\item [(b)] See (b) in Prop 2.3 in~\cite{Mari1}. 
	\end{description}
\end{proof}

\subsection{Admissible subgroups}
In this subsection we briefly recall the definition of admissible subgroups and the results given in~\cite{Mari2} for  the cyclic groups of order $p$ in $\mathrm{PGL}_{k+1}(\mathbb{C})$ for $p\geq k+1$
\begin{definition}
A finite subgroup $H$  of $\mathrm{PGL}_{k+1}(\mathbb{C})$ is called admissible if there exists a set of points  in  $\mathbb{P}^k(\mathbb{C})$ in general position and of cardinality $n\geq k+3$ such that the automorphism group of this point set contains $H$.
\end{definition}
We are interested in admissible cyclic subgroups. Recall that a finite cyclic subgroup of order $m$ of  $\mathrm{PGL}_{k+1}(\mathbb{C})$ can be put in diagonal form and  that any of its generator can be represented by a diagonal matrix of finite order $m$ whose entries are all powers of a primitive $m$-th root  
In~\cite{Mari2} we have shown that if the order of the group is a prime $p$ with $p\geq k+1$ a necessary and sufficient condition for  this groups to be  admissible is that the entries of the diagonal are all distinct. 
 	\section{Cyclic groups of order $p$ with
 $p\geq k+1$}\label{section1}

\subsection{Preliminary results}\label{auxiliary}

In order to  count the conjugacy classes of admissible cyclic subgroups of  order $p$ in $\mathrm{PGL}_{k+1}(\mathbb{C})$ when $p\geq k+1$ we need to make some preliminary computations.
In Section~\ref{section0} we have seen that an admissible subgroup $C_p$ of  $\mathrm{PGL}_{k+1}(\mathbb{C})$ can be put in diagonal form and  that any generator is represented by a matrix of finite order $p$ whose entries are all powers of a primitive $p$-th root  of unity $\epsilon_p$ and are all distinct. Up to rescaling this matrix we may assume that the last entry is equal to  one. Now,  taking a suitable power of this matrix  we get a matrix   whose  second to last  entry is $\epsilon_p$. It is clear that any subgroup $C_p$ contains only one such generator and that it is uniquely identified by such a matrix. Hereafter in this paper  such a matrix  will  be referred to as the identifying matrix, or in short as the id. matrix  \\
Then, the first step for  counting the number of conjugacy classes of admissible subgroups of order $p$ is then to count the number of identifying matrices. This is equivalent to count certain $(k+1)$-tuples of integers.
Indeed, fixing a primitive $p$-th root of unity $\epsilon_p$ one gets an isomorphism 
\begin{equation*}
\begin{split}
\varphi: & \mu_p\rightarrow \mathbb{Z}_p\\
& e_p^\alpha\mapsto \left[\alpha \right] 
\end{split} 
\end{equation*}
where $\mu_p$ is the multiplicative subgroups of $\mathbb{C}^*$ consisting of $p$-th roots of unity and $\mathbb{Z}_p$ is the additive group of the integers modulo $p$.
This isomorphism induces in turn an isomorphism
\begin{equation*}
\begin{split}
\varphi_{k+1} : E_{k+1}^p\rightarrow \mathbb{Z}_p^{\oplus (k+1)} 
\end{split} 
\end{equation*}
between the subgroup of $E_{k+1}^p$ of diagonal matrices  in $\mathrm{GL}_{k+1}(\mathbb{C})$ of  order dividing $p$  and the direct sum of $k+1$ copies of $\mathbb{Z}_p$.
In addition since any automorphism of $\mathrm{PGL}_{k+1}(\mathbb{C})$ of order $p$ in diagonal form is represented by a diagonal matrix of order $p$ and any such representations is given up to a multiplication by a $p$-th root of unity, we can identify these automorphisms with the quotient group $E_{k+1}^p/\bar{\mu}_p$, where  $\bar{\mu}_p$ are the scalar matrices which are the product of the identity matrix by a $p$-th  root of unity. We have also an automorphism between the quotient groups 
\begin{equation*}
\begin{split}
\bar{\varphi}_{k+1} : E_{k+1}^p/\bar{\mu}_p\rightarrow \mathbb{Z}_p^{\oplus (k+1)}/\langle I\rangle
\end{split} 
\end{equation*}
where $\langle I\rangle$ is the subgroup of  $\mathbb{Z}_p^{\oplus (k+1)}$ generated by the element whose entries are all equal to $1$. Let us denote the quotient group $\mathbb{Z}_p^{\oplus (k+1)}/\langle I\rangle$ with $X$.  By the assumptions made at the beginning of this section we can identify any generator of an admissible subgroup of order $p$ in 
$\mathrm{PGL}_{k+1}(\mathbb{C})$ with  an element $\bar{a}$ of $X$, whose representatives in $\mathbb{Z}_p^{\oplus (k+1)}$ have distinct entries.  In addition up to replacing $\bar{a}$ with a suitable  multiple, that is with a another generator of the additive subgroups $A$ generated by  $\bar{a}$ in $X$, we can choose as representative  for $a$ (and hence for $A$)  the only representative  whose last  entry is equal  to  $0$ and the second to the last is equal to $1$.
With an abuse of language,  we will refer to the tuples of these shape as identifying vectors, or in in short as id. vectors. With a slightly abuse of notation we will denote an identifying  vector with \begin{equation}\label{eqn:blockmatrix1}
\left[
\begin{array}{c}
a_1 \\ 
a_2 \\ 
\vdots \\
a_{k-1} \\
1 \\
0\\
\end{array}\right]
\end{equation}
or indifferently with $[{a_1,a_2, \cdots a_{k-1},1,0}]^T$, where the $a_i$ are distinct integers taken in the set $\lbrace{2, \ldots, p-2}\rbrace$. The same notation, without  the restriction on the entries for the identifying vector, will be used for the other elements in the quotient group $X$.\\
Going back to the problem of counting the conjugacy classes of finite admissible subgroups of order $p$ of  $\mathrm{PGL}_{k+1}(\mathbb{C})$, the above discussion shows  that the first step in this analysis, that is to count the number of identifying matrices, can be reduced to counting the id. vectors. Now, it is easy to see that this number
is the number of $k-1$-permutations of $p-2$ objects.  We indicate this number with $P(p-2,k-1)$.\\
Since we are dealing with conjugacy classes, the second step in our discussion is to consider the action of  the symmetric group  $S_{k+1}$ on the set of subgroups of $\mathrm{PGL}_{k+1}(\mathbb{C})$  of order $p$ and in diagonal form. Recall that the symmetric group $S_{k+1}$ can be identified as the group of permutation matrices in $\mathrm{PGL}_{k+1}(\mathbb{C})$. If  $H$ is the stabilizer  of an admissible group $G$ for the action of $S_{k+1}$ then the number of its diagonal conjugates equals the index $\left[S_{k+1}:H\right]$ of $H$ in $S_{k+1}$. By the identification subgroups with id. vectors,  this action of $S_{k+1}$ can be better explained in terms of these last ones. 
Let $a$ be an id. vector, $\bar{a}$ its class in $X$ and $A$ the additive group generated by this class. If the stabilizer $H$ in $S_{k+1}$ of $A$ has order $h$ there are $h$ classes in the additive group $A$ that are represented by  conjugates of $a$.  Moreover, the index $\left|  S_{k+1}:H\right|$ of $H$ in $S_{k+1}$ indicates the number of distinct conjugated additive subgroups which are generated by the classes of conjugates of $a$. Hereafter, these groups will be referred to as groups as conjugated subgroups of $A$ and the whole set as the conjugacy class of $A$. Note, also, that each of these groups has  an id. vector and contain  $h$  elements which are represented by conjugates of $a$.  To simplify the notation we will denote the index  $\left|  S_{k+1}:H\right|$ simply with $m_H$.

\begin{example}\label{example1}
	Let $k=2$ and $p=5$. Consider the  additive subgroup $A$ of $X$ generated by the class of the identifying  vector
	\begin{equation}\label{eqn:blockmatrix2}
	a=\left[
	\begin{array}{c}
	2 \\ 
	1 \\ 
	0\\
	\end{array}\right]
	\end{equation}
	This group consists of other $4$ classes which are represented by the following $3$-tuples
	\begin{equation}\label{eqn:blockmatrix3}
	\left[
	\begin{array}{c}
	4 \\ 
	2 \\ 
	0\\
	\end{array}\right], \ 
	\left[
	\begin{array}{c}
	1 \\ 
	3 \\ 
	0\\
	\end{array}\right], \ 
	\left[
	\begin{array}{c}
	3 \\ 
	4 \\ 
	0\\
	\end{array}\right], \  \left[
	\begin{array}{c}
	0 \\ 
	0 \\ 
	0\\
	\end{array}\right]
	\end{equation}
Since the pairs $\left[{3,4,0}\right]^T$ and $\left[{0,1,2}\right]^T$ and $\left[{1,3,0}\right]^T$ and $\left[{0,2,4}\right]^T$ lie respectively in the same class, the stabilizer of $A$ in $S_3$ contains  the group $H=\langle(13)\rangle$. A brief inspection shows that no other permutations in $S_3$ stabilize $A$, so the stabilizer is indeed $H$. Since $\left|  S_3:H\right|=m_H=3$,  we have $3$ lateral classes, respectively $H$, $(12)H$ and $(123)H$. Hence  the conjugates of $a$ split correspondingly in three sets
\begin{itemize}
	\item $\left[{2,1,0}\right]^T$,  $\left[{0,1,2}\right]^T$
	\item $\left[{1,2,0}\right]^T$,  $\left[{1,0,2}\right]^T$
	\item $\left[{0,2,1}\right]^T$,  $\left[{2,0,1}\right]^T$.
\end{itemize}
which corresponds to three distinct subgroups in $X$. Note in addition that taking the third multiple of  $\left[{1,2,0}\right]^T$ one get the id. vector of the second group which is $\left[{3,1,0}\right]^T$. Similarly, rescaling  $\left[{0,2,1}\right]^T$ one gets $\left[{4,1,0}\right]^T$ the id. vector of the third group.

\end{example}

The index $m_H$ helps to count the conjugacy classes, but before other steps are needed. The first is to notice that any permutation acting on the entries of each identifying vector determines a system of congruences modulo $p$, when considering these entries as unknowns.  The solutions of these congruences are the values for the entries $\alpha_1,\ldots, \alpha_{k-1}$ such that stabilizer in $S_{k+1}$   of the group generated by the  vector contains the cyclic subgroup  of $S_{k+1}$ generated by this permutation. Since we have imposed some conditions on the entries of the id. vector at the beginning of this section, only the solutions satisfying the criteria for the entries have to be kept. The solutions with this property will be called  admissible solutions. As we will see in some cases a system may admit solutions but no admissible solutions.
\begin{example}
Let $k=3$ and $p=7$. Let \begin{equation}\label{eqn:blockmatrix4}
	\left[
	\begin{array}{c}
	s \\ 
	1 \\ 
	0\\
	\end{array}\right]
	\end{equation}
be an identifying vector whose first entry is an unknown $s$. We want to find  for what values of $s$ the group  generated by the class of this vector  is stabilized  by the permutation $\alpha=(12)$ in $S_3$. Applying this permutation to the identifying vector one gets 
\begin{equation}\label{eqn:blockmatrix5}
\left[
\begin{array}{c}
1 \\ 
s \\ 
0\\
\end{array}\right]
\end{equation}
and this leads to the congruence  $s^2\equiv 1  \ (mod \ 7)$,  which has two solutions $s\equiv  1 \ (mod \ 7)$ and $s\equiv  -1 \ (mod \ 7)$. By our assumptions on the solutions only the second solution is admissible. 
\end{example}
Note also that everything here is up to  conjugation, that is if a system of congruences for a subgroup  in $S_{k+1}$ has admissible solutions so has any system relative to any subgroup in the conjugacy class of the original subgroup and the number of solutions is the same. 
It is clear that that for non-cyclic  subgroup of $S_{k+1}$,  a solution for  the  system of congruences relative to this subgroup is a solution common between the congruences of all generators.
 This last remark introduces what is left out in our discussion so  far, that is the  possibility of solutions in common between  systems of congruences for distinct subgroups of $S_{k+1}$.    So, when dealing with the solutions for  a subgroups one has to take out those in common with another subgroup.  Hypothetically one should consider all the ascending chains of subgroups in $S_{k+1}$ and starting from the lower group in the chain different from the identical one to see if there is any overlapping of solutions with the groups which stand higher in the chain. In reality our work here is greatly simplified. The first reason is that everything here is up to conjugation. Secondly our assumption on the solutions, i.e. we look for admissible solutions, will drastically reduce the number of conjugacy classes of subgroups involved. The classes retained will be referred to as the admissible classes. Anyway when in presence of solutions in common one has to remove these solution from the solutions in the smaller group.  \\
Completed this operation one has to multiply the number of solutions for an admissible subgroup  $H$ in $S_{k+1}$ by the number of subgroups in the conjugacy class of $H$: let us denote the number so obtained as $n_{H}$. Then, one has to compute the index $\left|  S_{k+1}:H\right|= m_H$ and  divide  $n_{H}$ by  $m_H$ to obtain the number of conjugacy classes of additive subgroups of $X$ stabilized by the group $H$. In terms of admissible subgroups of $\mathrm{PGL}_{k+1}\left(\mathbb{C}\right)$ of order $p$, this is the number of conjugate classes of these subgroups whose normalizer contains a subgroup in the conjugacy class of $H$ (see also Example~\ref{explicitcomputation} for a specific case showing computations)
What it has been done for $H$ has to be done for a representative of any admissible conjugacy class.\\
Finally, subtracting from $P(p-2,k-1)$, the total number of id. vectors,  the sum of all $n_H$ for $H$ ranging on all admissible conjugacy class and dividing the number so obtained by $k+1!$ one obtains the number of conjugacy classes of admissible subgroups of $\mathrm{PGL}_{k+1}\left(\mathbb{C}\right)$ of order $p$ whose normalizer contains no permutation. \\
Taking everything into account what we have seen so  far,  our  approach to the problem would proceed through the following steps:
\begin{description}[leftmargin=1.5 pc]
	\item[$(1)$]Determine and solve the systems of congruences relative to a representative of each conjugacy class of the $l$-cycles for $2\leq l\leq k+1$ and of some other suitable permutations which are the product of disjoint cycles.
	\item[$(2)$] Inspect the full list of conjugacy classes of $S_{k+1}$. Then on the basis  of the results in $(1)$ see which of these classes have to be rejected.
	\item[$(3)$] For each admissible conjugacy class of subgroups of $S_{k+1}$ determine the total number of solutions. Eventually neglect the solutions which are in common with a larger group. Finally, divide the total so obtained  by  $m_H=\frac{k+1!}{h}$, where $h$ is the cardinality of any representative of the class in question.  This  gives for each conjugacy class of proper subgroups of $S_{k+1}$ the  number of the corresponding conjugacy classes of admissible cyclic subgroups of order $p$ in $\mathrm{PGL}_{k+1}\left(\mathbb{C}\right)$ whose normalizer contains a subgroup in the conjugacy class of $H$.
	\item[$(5)$] Finally subtract the whole number of solutions of all the conjugacy classes of proper subgroups from  $P(p-2,k-1)$ and divide this number by $(k+1)!$. This gives the number of conjugacy classes whose normalizer contains no permutation.
\end{description}

\begin{remark}\label{terminology}
In the rest of this paper we will use the terminology and notation introduced in this section.  Moreover throughout this paper we will refer to some basic notions in the Theory of Groups  without introducing  or defining them. Examples  of these notions are the direct and the semidirect products of groups as well the definition of alternating or symmetric group. The reader can  refer to~\cite{Rob} for  the  full description of these concepts. 
\end{remark}
 
\begin{remark}
As we will see the correspondence between groups,  identifying matrices and identifying vectors, with the necessary adjustments, still holds when one replaces $p\geq k+1$ with $p< k+1$.
\end{remark}
\subsection{$k=2$}\label{case2}
The starting point of our analysis is the description of the finite imprimitive subgroups of $\mathrm{PGL}_{3}(\mathbb{C})$  given by Dolgachev and  Iskovskikh in~\cite{Dolgar}. In that paper the two authors consider a finite cyclic  subgroup $C_n$  of $\mathrm{PGL}_{3}(\mathbb{C})$  of order  $n\geq 3$ whose generator is represented by a matrix of finite order $n$, namely
\begin{equation}\label{eqn:blockmatrix6}
\left(
\begin{array}{c c c }
\varepsilon_n^s & 0 & 0  \\ 
0 & \varepsilon_n  & 0  \\ 
0	 & 0  & 1 \\ 
\end{array}\right)
\end{equation}
where $\varepsilon_n$ is a primitive $n$-th root of unity and $s>1$ is a positive integer. Then, they consider the action of the symmetric group $S_3$ on this subgroup. Recall that the symmetric group $S_3$ can be represented in  $\mathrm{PGL}_{3}(\mathbb{C})$ by the permutation matrices. The aim of their research is to find conditions on $n$ and $s$ such that the normalizer of $C_n$ in  $\mathrm{PGL}_{3}(\mathbb{C})$ contains some non-trivial subgroup of $S_3$. The outcome is that the normalizer of $C_n$ contains the permutation matrices of order $3$ if $s^2-s+1\equiv 0 \ (mod \ n)$ and the permutation $(12)$ if $s^2\equiv 1 \ (mod \ n)$. The two authors for their purposes do not need to find  the conditions for the remaining two transpositions in $S_3$ , but a straightforward computation shows that these are $s\equiv 2 \ (mod \ n)$ for transposition $(13)$ and $2s\equiv 1 \ (mod \ n)$ for transposition $(23)$.
Now, our aim is  to count the conjugacy classes of admissible finite cyclic subgroups of $\mathrm{PGL}_{3}\left(\mathbb{C}\right)$ of order $p$ with $p\geq 3$ corresponding to an abstract cyclic group of  order $p$. In Section~\ref{section0} we have shown that the admissible subgroups have generator represented by a diagonal matrix of finite order $p$ whose entries are all distinct and are  powers  of a primitive $p$-th root of unity. For our purposes we can apply what found out by the authors in~\cite{Dolgar} just replacing  $n$ with a prime $p\geq 3$ and explicitly assuming  that $s\not\equiv 0, 1 \ (mod \ p)$.  
We enlist in the following table all the above mentioned congruences for the permutations in  $S_3$.
\begin{table}[H]\label{table1}
	\centering
	\caption{Permutations in $S_3$}
	\begin{tabular}{|c|c|}
		\hline 
		Cycle & Congruence \\ 
		\hline 
		$(123)$  & $s^2-s+1\equiv 0 \ (mod \ p)$ \\ 
		\hline 
		$(12)$ & $s\equiv -1 \ (mod \ p)$ \\ 
		\hline 
		$(13)$& $s\equiv 2 \ (mod \ p) $\\ 
		\hline 
		$(23)$ & $2s\equiv 1 \ (mod \ p)$  \\ 
		\hline 
	\end{tabular} 
\end{table}
Analyzing the data in Table~\ref{table1}  one sees that when $p>3$ the congruence $s^2-s+1\equiv 0 \ (mod \ p)$ admits solutions only if $p\equiv 1 \ (mod \ 3)$ and that in this  case there are two distinct solutions. On the other hand, if $p=3$ then $s^3+1\equiv (s+1)^3 \ (mod \ 3)$ and the congruence $s^2-s+1\equiv 0 \ (mod \ p)$ admits a unique solution which is $s=2$. Note that  when $p=3$ $s=2$ is the common solution to of all the congruences in Table~\ref{table1}. A direct inspection also shows that apart for  the case $p=3$, for all other values of $p\geq 3$ there is no value of $s$ satisfying simultaneously more than one congruence in Table~\ref{table1}. 

\begin{theorem}\label{1theorem}
Let $p>3$ be a prime. Then the number of conjugacy classes of admissible subgroups of order $p$ in $\mathrm{PGL}_{3}\left(\mathbb{C}\right)$ is as follows depending on the parity of $p$ modulo $3$:
\begin{table}[H]
\centering
\label{table2}
	\caption{Formulae for conjugacy classes in $\mathrm{PGL}_{3}\left(\mathbb{C}\right)$}
\renewcommand{\arraystretch}{2}
\begin{tabular}{|c|c|c|}
	
	\hline 
Case &	Prime Parity &  Number of Conjugacy Classes \\ 
	\hline 
$(a)$ &  $p\not\equiv 1 \ (mod \ 3)$ & $1+ \frac{P(p-2,1)-3}{6}$ \\
	\hline 
$(b)$  & $p\equiv 1 \ (mod \ 3)$ &  $2+  \frac{P(p-2,1)-2-3}{6}$ \\ 
	\hline 
\end{tabular} 
\end{table}
\end{theorem}
\begin{proof}
\begin{description}
\item[(a)] Since  $p\not\equiv 1 \ (mod \ 3)$ the only congruences which admit solutions are those for the $2$-cycles. The subgroup $H$ generated by the $2$-cycle $\alpha=(12)$ is cyclic of order $2$. Its conjugacy class has length $3$. Since we have just one solution for each subgroup in this class and the solutions are all distinct we have $n_H=3$.  Dividing this number by $m_H=3$, we see that there just conjugacy class of admissible subgroups of order $p$ in $\mathrm{PGL}_{3}\left(\mathbb{C}\right)$ whose normalizer contains a transposition of $S_3$. For the remaining conjugacy classes it is enough to subtract $3$ from $P(p-2,1)$ and divide the total by $6$.
\item[(b)] Since  $p\equiv 1 \ (mod \ 3)$ all the congruences in Table~\ref{table1} have solutions. In particular the congruence for the $3$-cycle $(123)$ has $2$ distinct solutions. The conjugacy class of $H=\langle(123)\rangle$ in $S_3$ has length $1$, so we have $n_H=2$. Since $m_H=2$ we have just one conjugacy of subgroups of order $p$ in $\mathrm{PGL}_{3}\left(\mathbb{C}\right)$ class whose normalizer contains $H$. Since, the case of the $2$-cycles has already been discussed in $(a)$ above, to complete the proof it is enough to count the number of conjugacy classes whose normalizer contains no permutation: this number is given by subtracting the solutions for $H$ and for the $2$-cycles from $P(p-2,1)$ and divide this number by $6$.
\end{description}
\end{proof}
\begin{example}
	For each of the cases in Table~\ref{table2} we give some examples of explicit computation of the number of conjugacy classes. We also list for each example the number of elements in each class.
	\begin{table}[H]
		\centering
		\label{table3}
		\caption{Example}
		\begin{tabular}{|c|c|c|}
			\hline 
			Prime  &   Conjugacy Classes & $m_H$ \\
			\hline
			$p=5$ &  	 $1$ & $3$	\\
			\hline
			\multirow{2}{*}{$p=7$} &  $1$ & $2$ \\ &  $1$ & $3$ \\
			\hline  
			\multirow{2}{*}{$p=17$} & 	 $1$ & $3$	 \\  & 	 $2$ & $6$\\  
			\hline 
			\multirow{3}{*}{$p=19$} &  $1$ & $2$ \\ &  $1$ & $3$  \\
			&  $2$ & $6$  \\
			\hline  
		\end{tabular} 
	\end{table}
	
\end{example}
 \subsection{$k=3$}\label{case3}
According to the convention established in Remark~\ref{terminology} we describe a generator of an admissible  subgroup $C_p$  of order $p$  in $\mathrm{PGL}_{4}\left(\mathbb{C}\right)$ as an id. vector
\begin{equation}\label{eqn:blockmatrix7}
\left[
\begin{array}{c}
	s \\ 
	 t  \\ 
	1 \\ 
	0\\
\end{array}\right]
\end{equation}
where $s$ and $t$ are unknowns.
We consider initially  the $l$-cycles of $S_4$ for $2\leq l\leq 4$. As we have already seen it is  possible to restrict our attention to  a representative in each conjugacy class  of these cycles. 
 Under the assumption $s, t \not\equiv 0,1 \ (mod \ p)$, simple computations give the following conditions for each of the above mentioned cycles
\begin{table}[H]
	\centering
		\caption{Cycles in $S_4$}
		\label{table4}
	\begin{tabular}{|c|c|}
		\hline 
		Cycle & Congruence \\ 
		\hline 
		$(12)$ & $s\equiv t \ (mod \ p)$ \\ 
		\hline 
		$(123)$ & $t^2+t+1\equiv 0  \ (mod \ p)$, \ $s\equiv t^2  \ (mod \ p)$  \\
		 		\hline 
		$(1234)$ & $s^2\equiv -1 \ (mod \ p)$, \ $t\equiv s +1\ (mod \ p)$    \\ 
		\hline 
	\end{tabular} 
\end{table}

Analyzing the results in the above table, one sees that since we are supposing $s\not\equiv t \ (mod \ p)$ the transposition $(12)$ has to be rejected.  With regard to the $3$-cycle  	$(123)$, the relative system of congruences has solution only if $p\equiv 1\ (mod \ 3)$, as the first congruence requires $t$ to be a primitive $3$-rd root of unity modulo $p$.     Let us inspect now the results for the $4$-cycle $(1234)$. Here it is clear that the relative system of congruences has solution only if $p\equiv 1\ (mod \ 4)$ since the first congruence tells that $s$ has to be a primitive $4$-th root of unity modulo $p$.  \\
Then we consider the $3$ permutations in $S_4$ which are the product of two disjoint $2$-cycles. Though they lie in a single conjugacy class, for our purposes  we need to consider all  $3$ here. The corresponding congruences are listed in the following table

\begin{table}[h]
\centering
\caption{$2\times 2$ permutations in $S_4$}
\label{table5}
\begin{tabular}{|c|c|}
	\hline 
	Permutation & Congruence \\ 
	\hline 
	$(12) (34)$ & $t\equiv 1-s \ (mod \ p)$ \\ 
	\hline 
	$(13) (24)$ &  $t\equiv s+1  \ (mod \ p)$  \\
	\hline 
	$(14) (23)$,  & $t\equiv s-1 \ (mod \ p)$   \\ 
	\hline 
\end{tabular} 

\end{table}
Inspecting the congruences in Table~\ref{table5} it is clear that the congruence corresponding to the permutation  $(13) (24)$  is one of the two congruences for the $4$-cycle in Table~\ref{table4}, as the second power of  $(1234)$ is exactly $(13) (24)$.
Moreover a simple computation shows that for $p\geq 5$ no pair of values $(s,t)$ satisfying one of the congruences in Table~\ref{table5} satisfies any of the other two. 
For our purposes we have now to see  if there is any overlapping between the solutions of the congruences in Table~\ref{table4} and in Table~\ref{table5} and more in general if for certain values of $s$ and $t$ the normalizer of $C_p$  contains subgroups other than cyclic subgroups  generated by the $4$-cycles, $3$-cycles and the permutations which are the product of two disjoint $2$-cycles. In order to do this we first recall the list of conjugacy classes of proper subgroups of $S_4$. Note that in Table~\ref{table6} the first column indicates the isomorphism class  while the second stands for the identifier of the conjugacy class. If there is just one conjugacy class for one isomorphism class we keep the identifier used for  the isomorphism class.  We also sometimes with abuse of terminology and notation identify a representative of a conjugacy class with the identifier of the class. Note also that in Table~\ref{table6} we have neglected the conjugacy class of the trivial subgroup. For more reference on the set of conjugacy classes of subgroups of $S_4$ see~\cite{Sulai} or the list provided by the software G.A.P.(see~\cite{GAP4}).
\begin{table}[H]
	\centering
	\caption{Conjugacy classes in $S_4$}
	\label{table6}
	\renewcommand{\arraystretch}{2}
	\begin{tabular}{|c|c|c|c|c|}
	\hline 
   Group & Id. & Classes  &   Repres. & Length  \\ 
	\hline 
	\multirow{2}{*}{$C_2$} 	& $C_2^1$	&  $1$ &  $\langle(12)\rangle$ & $6$ \\ 
	& $C_2^2$  & $1$ & $\langle(13)(24)\rangle$ & $3$ \\ 
	\hline 
	$C_3$ &   & $1$ &  $\langle(123)\rangle$ & $4$ \\
	\hline 
	$C_4 $ &  &   $1$ &  $\langle(1234)\rangle$ &  $3$\\ 
	\hline 
	\multirow{2}{*}{$V_4$} &  $V_4^1$	&  $1$ &  $\langle(12), (34)\rangle$ & $3$ \\
		&  $V_4^2$ 	&  $1$ &   $\langle(12)(34),(13)(24)\rangle$	& $1$\\
	\hline 
	
	$S_3$ & 	&   $1$ & $\langle(12),(123)\rangle$	& $4$ \\
	\hline 
	$D_8$ & 	&  $1$ & $\langle(13),(1234)\rangle$	& $3$ \\
	\hline 
	$A_4$ &	& $1$ &  $\langle(123),(12)(34)\rangle$ & $1$ \\
	\hline 
	\end{tabular} 
	\end{table}
Analyzing the list of representatives of  each conjugacy class in Table~\ref{table6} bearing in mind what discussed so far with regard to  the congruences in Table~\ref{table4} and Table~\ref{table5} ones sees what follows
\begin{description}
	\item [$(1)$] Since there are no values of $(s,t)$ for which the normalizer of $C_p$ can contain a $2$-cycle, the following conjugacy classes in $S_4$ have to be rejected: $S_4$, $C_2^1$, $V_4^1$, $S_3$ and $D_8$. 
	\item [$(2)$] Since there are no values of $(s,t)$ for which the normalizer of $C_p$ contains at the same time all the permutations which are the product of two $2$-cycles also the following subgroups of $S_4$ have to be rejected: $V_4^2$ and  $A_4$. Indeed, recall that $A_4$
	 contains $V_4^2$ and $V_4^2$ contains all three of the above mentioned permutations.
	\item [$(3)$] The facts in $(1)$ and $(2)$  above show in particular that if  $\alpha$ and $\beta$ are two   $4$-cycles (respectively two $3$-cycles)  in $S_4$ such that $\alpha\notin\langle\beta\rangle$ then  the relative systems of congruences have no  solution in common   as the group $\langle\alpha, \ \beta\rangle$ they would generate is the whole $S_4$ (respectively $A_4$). Similarly one sees that there is no possible overlapping of solutions between a permutation  which is the product of two $2$ disjoint $2$-cycles and any $4$-cycle or any $3$-cycle. 
	\end{description} 
Relying on what we have seen so far about the subgroups of $S_4$, we  count the number of conjugacy classes  of admissible subgroups of order $p$ in $\mathrm{PGL}_{4}\left(\mathbb{C}\right)$.
\begin{theorem}\label{acase3}
	Let $p\geq 5$ be a prime. The number of conjugacy classes of admissible cyclic subgroups of order $p$ in $\mathrm{PGL}_{4}(\mathbb{C})$ is as follows depending on the parity of $p$  modulo $4$ and $3$: 
\begin{table}[H]
	\centering
	\label{table7}
	\caption{Formulae for conjugacy classes in $\mathrm{PGL}_{4}\left(\mathbb{C}\right)$}
	\renewcommand{\arraystretch}{2}
	\begin{tabular}{|c|c|c|}
		\hline 
Case	&	Prime Parity &  Number of Conjugacy Classes \\ 
		\hline 
		$(a)$	&	$p\not\equiv 1 \ (mod \ 4)\wedge p\not\equiv 1 \ (mod \ 3)$ & $\frac{(p-3)\cdot 3}{12} +\frac{P(p-2,2)-(p-3)\cdot 3}{24}$  \\ 
		\hline 

$(b)$	&	$p\not\equiv 1 \ (mod \ 4)\wedge p\equiv 1 \ (mod \ 3)$  & $1+\frac{(p-3)\cdot 3}{12} +\frac{P(p-2,2)-8-(p-3)\cdot 3}{24} $\\ 
		\hline  
$(c)$	&	$p\equiv 1 \ (mod \ 4)\wedge p\not\equiv 1 \ (mod \ 3)$ & $1+\frac{(p-3-2)\cdot 3}{12} +\frac{P(p-2,2)-(p-3)\cdot 3}{24} $  \\ 
		\hline 
$(d)$	&	$p\equiv 1 \ (mod \ 4)\wedge p\equiv 1 \ (mod \ 3)$ & $2+\frac{(p-3-2)\cdot 3}{12} +\frac{P(p-2,2)-(p-3)\cdot 3-8}{24} $ \\ 
\hline
					\end{tabular} 
\end{table}
	\end{theorem}
\begin{proof}
\begin{description}
	\item[$(a)$] Since	$p\not\equiv 1 \ (mod \ 4)\wedge p\not\equiv 1 \ (mod \ 3)$ only the linear congruences have solutions. Taking for instance that relative to the permutation $(13)(24)$ one sees that it has $p-3$ solutions. So the total number of
solutions for these permutations is $3\cdot (p-3)$. Since the solutions are conjugated in sets of $12$, the total number of conjugacy classes of subgroups, whose normalizer contains a cyclic group of order $2$ is $\frac{3\cdot (p-3)}{12}$. It remains
to count the number of conjugacy classes whose normalizer contains no
permutations. It is easy to see that this number is $\frac{(p-3)\cdot 3}{12} +\frac{P(p-2,2)-(p-3)\cdot 3}{24}$ and
this completes the proof for $(a)$.
\item [$(b)$] The first part of the proof it the same as that of case $(a)$. Here as  $p\equiv 1 \ (mod \ 3)$ we have in addition the solutions of the congruences regarding the $3$-cycles. Since the length of the conjugacy class of the subgroup
$\langle(123)\rangle$ is $4$ and there are two solutions for each subgroup, the total
number of these solutions is $8$. Since the solutions are conjugated in sets of $8$, it is clear that there is just one conjugacy class. So, it is enough to add $1$ at the beginning of the formula 
for case $(a)$ and subtract  $8$ from the numerator of the last fraction in the
same formula.
\item [$(c)$] Since $p\equiv 1 \ (mod \ 4)$ the congruence relative to the subgroup $\langle (1234)\rangle$ , has $2$ distinct solutions. Moreover the conjugacy class of  $\langle(1234) \rangle$ has length $3$. So, we have a total of $6$ solutions. Since they are conjugated in group of $6$ we have just one conjugacy class. Now we need to count the solutions for the permutation which are the product of $2$ disjoint $2$-cycles, taking into account that $6$ of these solutions have already been counted when dealing with the $4$-cycles. Therefore the total number of solutions for these permutations  is $3\cdot (p-5)$. A simple computation then gives the last fraction in the formula.
\item [$(d)$] The proof for this case follows straightforward from those of the previous cases.

\end{description}	
\end{proof}	
\begin{example}
	For each of the cases listed in Theorem~\ref{acase3} we give an example of explicit computation of the number of the conjugacy classes. We also list for each example the number of elements in each class. Note that for $p=11$ and $p=13$ there is more than one conjugacy class of subgroups whose normalizer contains $C_2^2$.
	\begin{table}[H]
			\label{table8}
			\caption{Examples}
		\centering
	\begin{tabular}{|c|c|c|}
		
		\hline 
		Prime  &   Conjugacy Classes & $m_H$ \\
		\hline
		\multirow{2}{*}{$p=11$} &  $2$ & $12$ \\ &  $2$ & $24$  \\
			\hline 
				\multirow{2}{*}{$p=7$} &  $1$ & $8$ \\ &  $1$ & $12$ \\
				\hline 
				\multirow{3}{*}{$p=17$} & 	 $1$ & $6$	 \\  & 	 $3$ & $12$\\  & 	 $7$ & $24$\\
				\hline  
	\multirow{4}{*}{$p=13$}  &  	 $1$ & $6$	 \\  & $1$ & $8$\\ & $2$ & $12$\\  & $3$ &  $24$\\
				\hline  
				\end{tabular} 
\end{table}
	\end{example}
\normalsize
\begin{example}\label{explicitcomputation}
	On Table~\ref{table8} we have already listed the number of conjugacy classes for the case $p = 11$. Here we want to describe this case in more detail to further clarify with a meaningful example the action of the symmetric group on an id. vector and	the computation of the conjugacy classes. Using the formulas in Theorem~\ref{acase3} we have seen that for $p = 11$ there are two conjugacy with $m_H=24$ and two with $m_H=12$.
	Brief computations show that the classes of length $m_H=24$ are represented by the id.vectors	$v_1 = \left[2, 7, 1, 0\right] ^T$ and $v_2 = \left[ 9, 5, 1, 0\right]^T $, while the two classes  with  $m_H=12$ are represented	by the vectors $v_3 = \left[ 2, 3, 1, 0\right] ^T$ and $v_4 = \left[ 3, 4, 1, 0\right] ^T$. Using the results displayed on	Table~\ref{table5} we see that the vectors representing these two groups are closed under the	action of the permutation $(13)(24)$. Now, in $(a)$ of Theorem~\ref{acase3} we have shown that the congruence relative to this permutation has $p-3 = 11-3 = 8$ solutions. Since the conjugacy class of the group $H =< (13)(24) > $ has length $3$ (see~\ref{table6}), the total
	number of solutions for the permutations that are the product of $2$ distinct $2$-cycles is $n_H= 24$. On the other hand $|S_4 : H| = m_H = 12$ and since $\frac{n_H}{m_H}=2$ this accounts for the $2$ classes of length $12$. This means that acting with the group $S_4$ on the vectors $v_3$
	and $v_4$, then rescaling and finding a suitable multiple (if necessary) so that the second to the last entry is equal to $1$, one gets $12$ id. vectors. These vectors represent groups	which are conjugated to the group generated by the original vector and are stabilized 	by one of the groups in the conjugacy class of $H$. We list these vectors here below	displayed in groups of $4$ according to their stabilizer.
	\tiny 
	\begin{center}
		$\left[ \begin{array}{c}
		3\\
		4\\
		1\\
		0
	\end{array}\right] $ $\left[ \begin{array}{c}
	4\\
	5\\
	1\\
	0
	\end{array}\right] $  $\left[ \begin{array}{c}
	7\\
	8\\
	1\\
	0
	\end{array}\right] $  $\left[ \begin{array}{c}
	8\\
	9\\
	1\\
	0
	\end{array}\right] $ \hspace{0,1 cm}
		$\left[ \begin{array}{c}
		4\\
		3\\
		1\\
		0
		\end{array}\right] $ $\left[ \begin{array}{c}
		5\\
		4\\
		1\\
		0
		\end{array}\right] $  $\left[ \begin{array}{c}
		8\\
		7\\
		1\\
		0
		\end{array}\right] $  $\left[ \begin{array}{c}
		9\\
		8\\
		1\\
		0
		\end{array}\right] $ \hspace{0.1 cm}
			$\left[ \begin{array}{c}
			5\\
			7\\
			1\\
			0
			\end{array}\right] $ $\left[ \begin{array}{c}
			7\\
			5\\
			1\\
			0
			\end{array}\right] $  $\left[ \begin{array}{c}
			3\\
			9\\
			1\\
			0
			\end{array}\right] $  $\left[ \begin{array}{c}
			9\\
			3\\
			1\\
			0
			\end{array}\right] $\\ \vspace{0.5 cm}
	$\left[ \begin{array}{c}
	2\\
	3\\
	1\\
	0
	\end{array}\right] $ $\left[ \begin{array}{c}
	5\\
	6\\
	1\\
	0
	\end{array}\right] $  $\left[ \begin{array}{c}
	6\\
	7\\
	1\\
	0
	\end{array}\right] $  $\left[ \begin{array}{c}
	9\\
	10\\
	1\\
	0
	\end{array}\right] $ \hspace{0,1 cm}
		$\left[ \begin{array}{c}
		3\\
		2\\
		1\\
		0
		\end{array}\right] $ $\left[ \begin{array}{c}
		6\\
		5\\
		1\\
		0
		\end{array}\right] $  $\left[ \begin{array}{c}
		7\\
		6\\
		1\\
		0
		\end{array}\right] $  $\left[ \begin{array}{c}
		10\\
		9\\
		1\\
		0
		\end{array}\right]$ \hspace{0,1 cm}	$\left[ \begin{array}{c}
			4\\
			8\\
			1\\
			0
			\end{array}\right] $ $\left[ \begin{array}{c}
			8\\
			4\\
			1\\
			0
			\end{array}\right] $  $\left[ \begin{array}{c}
			2\\
			10\\
			1\\
			0
			\end{array}\right] $  $\left[ \begin{array}{c}
			10\\
			2\\
			1\\
			0
			\end{array}\right] $
	
	\end{center}
\end{example}

\normalsize 
\newpage
\subsection{$k=4$}
As usual we describe a generator of an admissible  subgroup $C_p$  of order $p$  in $\mathrm{PGL}_{5}\left(\mathbb{C}\right)$
as an id. vector
\begin{equation}\label{eqn:blockmatrix8}
\left[
\begin{array}{c}
s \\ 
t  \\ 
u\\
1 \\ 
0\\
\end{array}\right]
\end{equation}
where $s$, $t$ and $u$ unknowns. Here again in our computation we restrict our attention to a representative  for the each conjugacy class of  the $l$-cycles of $S_5$ for $2\leq l\leq 5$. For our purposes in our analysis we consider also some other suitable permutation.  Under the assumption $s, t, u \not\equiv 0,1 \ (mod \ p)$, straightforward computations give the following conditions for each of the above mentioned cycles
\begin{table}[H]
	\centering
	\caption{Permutations in $S_5$}
	\label{table9}
	\begin{tabular}{|c|c|}
		\hline 
		Cycle & Congruence \\ 
		\hline 
		$(12)$ & $s\equiv t \ (mod \ p)$ \\ 
		\hline 
		$(123)$ & $s\equiv t\equiv u  \ (mod \ p)$ \\
		\hline 
		\multirow{3}{*}{$(1243)$}  & $t^2+1\equiv 0 \ (mod \ p)$\\  & $s\equiv t^2\ (mod \ p)$\\  &  $u\equiv t^3 \ (mod \ p)$  \\ 
		\hline 
		\multirow{3}{*}{$(12345)$}	 & $s^4-s^3+s^2-s+1\equiv 0 \ (mod \ p)$\\
		& $t\equiv  -s^2 + s \ (mod \ p)$\\
		  &  $u\equiv s^3-s^2+s \ (mod \ p)$  \\ 
		\hline 
			\multirow{2}{*}{$(14)(23)$} 
		 & $s\equiv - 1 \ (mod \ p)$ \\
		 & $u\equiv -t \ (mod \ p)$ \\
			\hline
				\multirow{2}{*}{$(12)(34)$} 
				& $u\equiv - 1 \ (mod \ p)$ \\
				& $t\equiv -s \ (mod \ p)$ \\
				\hline
	\end{tabular} 
\end{table}

Inspecting the congruences in Table~\ref{table9}, one sees that since we are supposing $s\not\equiv t\not\equiv u \ (mod \ p)$ the cycles  $(12)$ and  $(123)$ have no solutions.  With regard to the $4$-cycle $(1243)$ one sees that the relative system of congruences has solution  only if $p\equiv 1\ (mod \ 4)$, as the first congruence requires $t$ to be a primitive $4$-th root of unity modulo $p$.  Note also that the value of $s$ is fixed and equal to $-1$, while $u$, being the third power of $t$ is a $4$-th primitive root of unity. With regard to the system of congruences relative to the $5$-cycle $(12345)$ one sees that it ha solutions  only if  $p\equiv 1\ (mod \ 5)$ as the congruence involving $s$ require $s$ to be a primitive $10$-th root of unity. Finally, please note that  both of the systems for the permutations which are the product of two $2$-cycles admit solutions and these solutions contain a free variable.\\
Now, we consider the list of conjugacy classes of proper subgroups of $S_5$ to see which of these classes can be neglected. Note also that in Table~\ref{table10} we have already excluded the conjugacy class of the trivial subgroup. For more reference on the set of conjugacy classes of subgroup of $S_5$ see the list found in~\cite{}  by the software G.A.P.(see~\cite{GAP4}) 
\begin{table}[H]
	\centering
	\caption{Conjugacy classes in $S_5$}
	\label{table10}
	\renewcommand{\arraystretch}{1.5}
	\begin{tabular}{|c|c|c|c|c|}
		\hline 
		Group & Id. & Classes  &   Repres. & Length  \\ 
		\hline 
		\multirow{2}{*}{$C_2$} 	& $C_2^1$	&  $1$ &  $\langle(12)\rangle$ & $10$ \\ 
		&  $C_2^2$  & $1$ & $\langle(14)(23)\rangle$ & $15$ \\ 
		\hline 
		$C_3$ &   & $1$ &  $\langle(123)\rangle$ & $10$ \\
		\hline 
		$C_4 $ &  &   $1$ &  $\langle(1243)\rangle$ &  $15$\\ 
		\hline 
		\multirow{2}{*}{$V_4$} &  $V_4^1$	&  $1$ &  $\langle(12), (34)\rangle$ & $15$ \\
		&  $V_4^2$ 	&  $1$ &   $\langle(14)(23),(13)(24)\rangle$	& $5$\\
		\hline 
		$C_5$ &  &   $1$ &  $\langle(12345)\rangle$ &  $6$\\ 
		\hline 
		$C_6 $ &  &   $1$ &  $\langle(12)(345)\rangle$ &  $10$\\
		\hline 
		\multirow{2}{*}{$S_3$} &  $S_3^1$	&  $1$ &  $\langle(12), (123)\rangle$ & $10$ \\
		& $S_3^2$	&   $1$ & $\langle(123),(12)(45)\rangle$	& $10$ \\
		\hline 
		$D_8$ &	&  $1$ & $\langle(13),(1234)\rangle$	& $15$ \\
		\hline 
		$D_{10}$ & 	&  $1$ & $\langle(12345),(14)(23)\rangle$	& $6$ \\
		\hline 
		$A_4$ & & $1$ &  $\langle(123),(12)(34)\rangle$ & $5$ \\
		\hline 
		$S_3\times S_2$ &	& $1$ &  $\langle(123),(12)(45)\rangle$ & $10$\\
		\hline 
		$D_{20}$ &	&  $1$ & $\langle(12345),(1243)\rangle$	& $6$ \\
		\hline 
		$S_4$ & 	&  $1$ & $\langle(12),(1234)\rangle$	& $5$ \\
		\hline 
		$A_5$ & 	&  $1$ & $\langle(12345),(13)\rangle$	& $1$ \\
		\hline 
		\end{tabular} 
		\end{table}

\begin{description}
\item [$(1)$] Since there are no values of  $(s,t)$  for which the normalizer of $C_p$ can contain
a $2$-cycle  the  following groups have be rejected: $C_2^1$, $V_4^1$, $C_6$, $S_3^1$, $D_8$, $S_4$ and $A_5$.
\item [$(2)$] Since there is no  solution for the congruences for the $3$-cycles also the following
subgroups of $S_5$ have to be excluded: $C_3$, $S_3^2$, $A_4$  and $S_3\times S_2$.
\item [$(3)$] Since the congruences for the two permutations $(14)(23)$ and $(12)(34)$ have no solution in common also the  group $V_4^2$ has to be rejected.
\item [$(4)$] With regard to the group $D_{10}$, one sees that in addition to the $5$-cycle $(12345)$ it contains also  the permutation $(14)(23)$. Since the congruences relative to  this permutation require $s\equiv -1\ (mod \ p)$   this value inserted in the congruence
$s^4-s^3+s^2-s+1\equiv 0 \ (mod \ p)$ forces $p=5$. A simple computation then shows that for this value of  $s$  the system of congruences for $(12345)$ has the solution $s\equiv 4 \ (mod \ 5)$, $t\equiv 3 \ (mod \ 5)$  and $u\equiv 3 \ (mod \ 5)$. Now,  this solution is also a solution for $(14)(23)$. Moreover a further inspection shows that this solutions satisfies also the system for $(1234)$ and hence the  congruences for both the generators of the representative of the conjugacy  class identified with $D_{20}$.

\end{description}
\begin{theorem}\label{case4}
Let $p>5$ be a prime .The number of conjugacy classes of admissible cyclic subgroups of order $p$ in $\mathrm{PGL}_{5}(\mathbb{C})$ is as follows depending on the parity of $p$ modulo  $5$ and $4$

\begin{table}[h]
	\centering
	\label{table11}
	\caption{ Formulae  for conjugacy classes in $\mathrm{PGL}_{5}\left(\mathbb{C}\right)$}
	\renewcommand{\arraystretch}{2}
	\begin{tabular}{|c|c|c|}
		\hline 
Case	&	Prime Parity &  Number of Conjugacy Classes \\ 
\hline	
$(a)$ & $p\not\equiv 1 \ (mod \ 5)\wedge p\not\equiv 1 \ (mod \ 4)$ & $\frac{(p-3)\cdot 15}{60} +\frac{P(p-2,3)-(p-3)\cdot 15}{120}$  \\ 
\hline 	
$(b)$ &	$p\not\equiv 1 \ (mod \ 5)\wedge p\equiv 1 \ (mod \ 4)$  & $1+\frac{(p-5)\cdot 15}{60} +\frac{P(p-2,3)-(p-3)\cdot 15}{120}$\\ 
\hline  
$(c)$& $p\equiv 1 \ (mod \ 5)\wedge p\not\equiv 1 \ (mod \ 4)$ & $1+\frac{(p-3)\cdot 15}{60} +\frac{P(p-2,3)-24-(p-3)\cdot 15}{120}$  \\ 
\hline 
$(d)$ & $p\equiv 1 \ (mod \ 5)\wedge p\equiv 1 \ (mod \ 4)$ & $2+\frac{(p-5)\cdot 15}{60} +\frac{P(p-2,3)-24+(p-3)\cdot 15}{120}$ \\ 
\hline 
\end{tabular} 
\end{table}
\end{theorem}
\begin{proof}
	\begin{description}
		\item[$(a)$] Since $p\not\equiv 1 \ (mod \ 5)\wedge p\not\equiv 1 \ (mod \ 4)$  the only systems  of congruences in Table~\ref{table10} which admit solutions are those for the permutation which are the product of $2$ distinct $2$-cycles. Taking for  instance that relative to the permutation $(12)(34)$ a simple computation shows that it has $p-3$ solutions. 
		So, the total number of solutions for these permutations is $15\cdot (p-3)$ as there $15$ of these permutations. Since these solutions are conjugated in sets of $60$, the total number of conjugacy classes of subgroups whose normalizer contains a cyclic group of order $2$ is $\frac{15\cdot (p-3)}{60}$. It remains to count the number of conjugacy classes whose normalizer contains no permutation. It is easy to see that this number is $\frac{P(p-2,3)-(p-3)\cdot 15}{120}$  and this completes the proof for $(a)$. 
\item[$(b)$] Since $p\equiv 1 \ (mod \ 4)$ the congruence relative to the subgroup $\langle(1234)\rangle$ has
$2$ distinct  solutions which are conjugated. Moreover the conjugacy class
of $\langle(1234)\rangle$ has length $15$. So, we have a total of $30$ solutions.
Since they ere conjugated in sets of $30$, there is exactly one conjugacy class. Now, we need to count the solutions for the permutations which
are the product of two disjoint $2$-cycles taking into account that $30$ of
these solutions have already been counted when dealing with the $4$-cycles.
Therefore their total number is $15\cdot (p-5)$. Then, the number of conjugacy classes whose normalizer contains no permutation is $\frac{P(p-2,3)-(p-3)\cdot 15}{120}$. Taking everything into account one gets the formula for this case $(b)$.
\item[$(c)$]  The first part of the proof is the same as that of case $(a)$ as $p\not\equiv 1 \ (mod \ 4)$.
Here in addition  we have  the solutions of the congruences relative to  the $5$-cycles.
Since the length of the  conjugacy class of the subgroup $\langle(12345)\rangle$ is $6$
and there are $4$ solutions for each subgroup the total number of solutions
is $24$. Since these solutions are conjugated in set of $24$, it is clear that
there in just one conjugacy class. So, it is enough to add  a $1$ in the formula
for  the case $(a)$ and subtract $24$ in the numerator of the last fraction in the
same formula.
\item[$(d)$] The proof for this case follows straightforward from those of the previous
cases. 
	\end{description}
\end{proof}
\begin{remark}
 It remains to consider the case $p=5$ . In this case $P(5-2,3)= 6$ and we
have already seen that there is a solution for the group $D_{20}$ whose conjugacy class
has length $6$. So,  there is just one conjugacy class.
\end{remark}

\begin{example}
	For each of the cases listed in Theorem~\ref{case4} we give an example of explicit computation of the number of conjugacy classes. We also list for each example the number of elements in each class.
	\begin{table}[H]
		\centering
		\begin{tabular}{|c|c|c|}
			\hline 
			Prime  &   Conjugacy Classes & $m_H$ \\
			\hline
			$p=7$ &  $1$ & $60$ \\
			\hline 
			\multirow{3}{*}{$p=13$} &  $1$ & $30$ \\ &  $2$ & $60$ \\ &  $7$ & $120$ \\
			\hline
			\multirow{3}{*}{$p=11$} & 	 $1$ & $24$	 \\  & 	 $2$ & $60$\\  & 	 $3$ & $120$\\
			\hline 
			\multirow{4}{*}{$p=41$}  &  	 $1$ & $24$	 \\  & $1$ & $30$\\ & $9$ & $60$\\  & $452$ &  $120$\\
			\hline
					\end{tabular} 
	\end{table}
\end{example}

\newpage
\subsection{$k=5$}

According to the convention introduced in Remark~\ref{terminology} we describe a generator of an admissible  subgroup $C_p$  of order $p$  in $\mathrm{PGL}_{6}\left(\mathbb{C}\right)$  as an id. vector

\begin{equation}\label{eqn:blockmatrix9}
\left[
\begin{array}{c}
s \\ 
t  \\ 
u\\
v\\
1 \\ 
0\\
\end{array}\right]
\end{equation}
where $s$, $t$, $u$ and $v$ are distinct integers with  $2 \leq s,t, u, v\leq p-1$.
 Here again for our purposes we consider a representative of each conjugacy class of $l$-cycles of $S_6$ for $2\leq l\leq 6$. We also  consider the representative of some other suitable class of permutations. Please note that in the following table $i\in\lbrace{\epsilon_3,\epsilon_3^2}\rbrace$ where $\epsilon_3$ is a primitive $3$-rd root of  unity modulo $p$.
\begin{table}[H]
	\centering
	\caption{Permutations}
	\label{table12}
	\begin{tabular}{|c|c|}
	\hline 
	Cycle & Congruence \\ 
	\hline 
	$(12)$ & $s\equiv t \ (mod \ p)$ \\ 
	\hline 
	$(123)$ & $s\equiv t\equiv u  \ (mod \ p)$ \\
	\hline 
	$(1234)$  & $s\equiv t\equiv u  \equiv v \ (mod \ p)$  \\ 
	\hline 
	\multirow{4}{*}{$(12345)$}	 & $v^4+v^3+v^2+v+1\equiv 0 \ (mod \ p)$\\
	& $s\equiv v^4 \ (mod \ p)$\\
	&  $u \equiv v^3 \ (mod \ p)$  \\ 
	&  $u\equiv v^2 \ (mod \ p)$\\
	\hline 
	\multirow{4}{*}{$(123456)$} 
	& $s^s+s+1\equiv 0 \ (mod \ p)$ \\
	& $t\equiv -s^2+s \ (mod \ p)$ \\
	& $u\equiv s^3-s^2+s \ (mod \ p)$ \\
	& $v\equiv s^4-s^3+s^2-s \ (mod \ p)$ \\
	\hline
		\multirow{3}{*}{$(135)(246)$} 
	& $s\equiv i^2- i^2\cdot v\ (mod \ p)$ \\
	& $t\equiv - i^2\cdot v \ (mod \ p)$ \\
		& $u\equiv i + v \ (mod \ p)$ \\
	\hline
	\multirow{2}{*}{$(14)(25)(36)$} 
	& $s\equiv u-v\ (mod \ p)$ \\
	& $t\equiv u-1 \ (mod \ p)$ \\
	\hline
	\end{tabular} 
	\end{table}
\begin{remark}\label{congruence5}	
The analysis of the congruences shown in Table~\ref{table12} leads to the following remarks:
\begin{description}
	\item[$(1)$] The congruences for the $2$-cycle, $3$-cycles and $4$-cycles as well as for the permutations which are the product of $2$ disjoint $2$-cycles have no solution under the requirement  $s\not\equiv t\not\equiv u \not\equiv v \ (mod \ p)$.
	This in turn shows that also the permutations which are the product of a $2$-cycle and a disjoint $3$-cycle or of a $2$-cycle and a disjoint $4$-cycle have to be rejected, as the second power of these permutations is respectively a $3$-cycle or a  permutation which is the product of $2$ disjoint $2$-cycles.
	\item[$(2)$]  The system of congruences for the permutation $(12345)$ has solution only if $p \equiv 1 \ (mod \ 5)$ has the first congruence in the system requires  $s$ to be a $5$-th primitive root of unity modulo $p$. Moreover, it is easy to see that the entries of each solution of the system are all distinct.
		\item[$(3)$]  The system of congruences for the cycle $(123456)$ has solution only if $p \equiv 1 \ (mod \ 3)$ as the first congruence in the system requires $s$ to be a $3$-rd primitive root of unity. Here again a simple inspection also shows that the entries of the solutions of the system are all distinct. 
	\item[$(4)$] The system of congruences for the permutation $(135)(246)$ admits solution only if $p \equiv 1 \ (mod \ 3)$ as the solution of each of congruences depends on a power of $\epsilon_3$, where $\epsilon_3$ stands for a primitive $3$-rd root of unity modulo $p$.  
\end{description} 
\end{remark}	

Unlike the  previous cases, we do not give here a detailed list of the conjugacy classes of $S_6$ as they are  $55$ in total. Another reason for this it is that many of these classes  can be straightforwardly rejected because their representatives  contain the permutations that have no solution described in $(1)$ of Remark~\ref{congruence5}. These classes include:
\begin{itemize}
\item  all  the classes of subgroups of order $2$ except those whose generator is a permutation which is the products of three $2$-cycles.
\item  the class of cyclic groups of order $3$ whose generator is a $3$-cycle.
\item all the classes of subgroups of order $4$, $8$, $9$, $10$, $12$, $16$, $18$, $20$, $24$ and $72$.
\item  all the classes  of subgroups of  order $6$ except two of these classes.
	\end{itemize}
For more reference on the set of conjugacy classes of subgroups of $S_5$ see~\cite{Vekri} or the list provided by the software G.A.P.(see~\cite{GAP4}).
Now, we need to investigate the two remaining classes of subgroups of order $6$. One is represented by cyclic groups of order $6$ generated by a $6$-cycle: this is the class represented by the $6$-cycle $(123456)$. So by the result in Table~\ref{table12} this  class is kept when $p \equiv 1 \ (mod \ 3)$. The second class of subgroups of order $6$ is that represented by the subgroup generated by the permutations $(135)(264)$ and $(14)(23)(56)$. This group is isomorphic to $S_3$ and contains also the permutation $(14)(25)(36)$ whom we have already dealt with in Table~\ref{table12}.  Now, a straightforward computation shows that the system of congruences for the permutation $(14)(23)(56)$  is 
\begin{equation*}\label{conditions}
\begin{cases}
s\equiv 1-v \ (mod \ p) \\
t\equiv 1-u \ (mod \ p).
\end{cases}
\end{equation*}
Comparing this system with that for $(14)(25)(36)$ one sees that if a value for $u$ satisfies both systems this value has to be $1$. So the class represented by this subgroup has to be rejected.
Taking everything into account we see that the only classes that we can retain are those represented by the following subgroups
\begin{table}[H]
	\centering
	\caption{Conjugacy classes}
	\label{table13}
	\begin{tabular}{|c|c|c|}
		\hline 
		Group & Repres. & Class. Length  \\ 
		\hline 
		$C_2$ 	  &  $\langle(14)(25)(36)\rangle$ & $15$ \\ 
		\hline 
		$C_3$ &  $\langle(135)(246)\rangle$ & $20$ \\
		\hline 
				$C_5 $ &  $\langle(12345)\rangle$ &  $36$\\ 
		\hline 
		$C_6 $ &  $\langle(123456)\rangle$ &  $120$\\
		\hline
		\end{tabular} 
\end{table}

This in turn shows that the only possible overlapping of solutions is between the solutions of the permutations belonging to the group generated by a $6$-cycle when $p \equiv 1 \ (mod \ 3)$.
For our purposes we need then to count the solutions for the system of congruences relative to the representatives of the admissible conjugacy classes identified with $C_2$ and $C_3$.

\begin{lemma}\label{counting}
	Let $p\geq 7$ be a prime.
	\begin{description}
	\item[$(a)$] The number of admissible solutions of the system of congruences  relative to	the permutation $(14)(25)(36)$ is $(p-3)\cdot (p-5)$.
		\item[$(b)$] If $p\equiv 1 \ (mod \ 3)$ the system of congruences for the permutation $(135)(246)$
	admits solutions  and the number of these solutions  is $2\cdot(p-4)$.
	\end{description}
\end{lemma}
\begin{proof}
	\begin{description}
		\item[$(a)$] By the results displayed in Table~\ref{table12} we look for solutions of the shape $[u-v,u-1, u, v,1,0]$.
	We first count the possible values for $u$. Since the entries of $[u-v,u-1, u, v,1,0]$ have to be all distinct, we have that $u$ has to satisfy simultaneously the conditions $2\leq u\leq p-1$ and $2\leq u-1 \leq p-1$. This shows that $3\leq u\leq p-1$ which gives a total of $p-3$ choices. Now we have to count the number of choices we have  for $v$. Clearly $v$  has to satisfy  $v\neq 0, 1, u, u-1$. Note that if $v$ satisfy these conditions so does $u-v$. So, we are left with the condition $u-v\neq v$ (or equivalently $2v\not\equiv u \ (mod \ p)$) It is easy to see that this condition is indipendent from the previous ones. Hence we have a total of $(p-5)$ choices for $v$. This leads to a total of $(p-3)\cdot (p-5)$ solutions for the system of congruences for $(14)(35)(26)$.
	\item [$(b)$] Since $p\equiv 1 \ (mod \ 3)$ the system of congruences for the permutation $(135)(246)$ admits solutions by the results  shown in Table~\ref{table12}. In order to count these solutions we have only to see for which values of $v\neq  0, 1$ the vector $[i^2- i^2\cdot v, - i^2\cdot v, i+v, v, 1, 0]$ (where $i\in\lbrace{\epsilon_3,\epsilon_3^2}\rbrace$) has distinct entries. Short computations shows that the requirements $i+v\not\equiv  0,1\ (mod \ p)$ cover all the needed conditions. These $2$ inequalities show that $v$ has to be distinct modulo $p$ from $-i$ and $1-i$. This leads to $(p-4)$ choices for $v$. Since there are  two choices for $i$ we have a total of $2\cdot (p-4)$ solutions.
\end{description}	
\end{proof}

\begin{theorem}\label{k5}
Let $p\geq 7$ be a prime. The number of conjugacy classes of admissible cyclic subgroups of order $p$ in $\mathrm{PGL}_{6}(\mathbb{C})$ is as follows depending on the parity of $p$ modulo  $5$ and $3$:
	\begin{table}[H]
	\label{table14}
	\caption{Formulae for the conjugacy classes in $\mathrm{PGL}_{6}\left(\mathbb{C}\right)$}
	\begin{center}
		\renewcommand{\arraystretch}{2}			
		\begin{tabular}{|c|c|c|} \hline
			Case	& Prime Parity & Number of Conjugacy Classes \\ \hline
				$(a)$ & $p\not\equiv 1 \ (mod \ 5)\wedge p\not\equiv 1 \ (mod \ 3)$ & $\frac{\alpha}{360} +\frac{P(p-2,4)-\alpha}{720}$\\ \hline
			$(b)$ & $p\equiv 1 \ (mod \ 5)\wedge p\not\equiv 1 \ (mod \ 3)$ & $1+\frac{\alpha}{360} +\frac{P(p-2,4)-144-\alpha}{720}$\\ \hline
						$(c)$ & $p\not\equiv 1 \ (mod \ 5)\wedge p\equiv 1 \ (mod \ 3)$ & $1+\frac{\beta}{360}+\frac{\delta}{240}+\frac{P(p-2,4)-120-\beta- \delta}{720}$
			\\ \hline
			$(d)$ & $p\equiv 1 \ (mod \ 5)\wedge p\equiv 1 \ (mod \ 3)$ & $2+\frac{\beta}{360}+\frac{\delta}{240} +\frac{P(p-2,4)-120-144- \beta-\delta }{720}$	\\ \hline 
		\multicolumn{3}{l} {where $\alpha=(p-3)\cdot(p-5)\cdot 15$, $\beta=(p-3)\cdot(p-5)\cdot 15-120$ and $\delta=2\cdot (p-4)\cdot 20-120$ }
			\end{tabular}
	\end{center}
\end{table}
\end{theorem}
\begin{proof}
\begin{description}
	\item[$(a)$] Since $p\not\equiv 1 \ (mod \ 5)\wedge p\not\equiv 1 \ (mod \ 3)$ in Table~\ref{table12} only the system  of congruences for the permutation which is the product of $3$ distinct $2$-cycles admit solutions. By $(a)$ in Lemma~\ref{counting} the total number of solutions for this system is $\alpha=(p-3)\cdot (p-5)$. Since the conjugacy class of this permutation contains $15$ elements, this gives a total of $15\cdot (p-3)\cdot (p-5)$ solutions. Since these solutions group in sets of $360$, the number of conjugacy classes of subgroups whose normalizer contains a cyclic group of order $2$
	  is  $\frac{(p-3)\cdot(p-5)\cdot 15}{360}$. It is then clear that the number of conjugacy classes whose normalizer contains no permutation is  $\frac{P(p-2,4)- \alpha}{720}$. So we are done with case $(a)$.
	  		\item[$(b)$] Since $p\equiv 1 \ (mod \ 5)\wedge p\not\equiv 1 \ (mod \ 3)$, in addition to the  solutions  for the product of $3$ distinct $2$-cycles  we have also the solutions for the subgroups generated by the $5$-cycles. We have $4$-solutions for each subgroup and there are $36$ subgroups in the class which gives a total of $144$ solutions.  Since these solutions go in groups of $144$, we have just $1$ conjugacy class of subgroups whose normalizer contains a cyclic group of order $5$. Relying on what we have seen in case $(a)$ above and taking into account these $144$ solutions, it is the straightforward to prove the rest of the formula for this case.
		\item[$(c)$] Since $p\not\equiv 1 \ (mod \ 5)\wedge p\equiv 1 \ (mod \ 3)$, three systems of congruences in  Table~\ref{table12} have solutions: that for the $6$-cycle, and those  for the permutations $(14)(25)(36)$ and  $(135)(246)$. The first system accounts for $2$ solutions.  Since the group generated by the $6$-cycle has $60$ conjugates this give $120$ solutions. These solutions group in sets of $120$, so we have $1$ conjugacy class. Now, in $(a)$ we have counted the solutions relative to the permutations in the conjugacy class of the group generated by $(14)(25)(36)$ and have seen that they are $15\cdot(p-3)\cdot (p-5)$. What changes here is that due to an overlapping solutions   we have to take $120$ solutions away from this number. Indeed, each of the $15$ subgroups  in the conjugacy class of $\langle(14)(25)(36)\rangle$ is contained in $4$ subgroups belonging to the conjugacy class of  $\langle(123456)\rangle$. These makes a total of $8$ solution for each of the $15$ subgroups which leads to a final number of $120$.  Subtracting $120$ from $15\cdot(p-3)\cdot (p-5)$ and dividing by $360$ we obtain that the number conjugacy classes of subgroups whose normalizer contains a cyclic group of order $2$
		is  $\frac{\beta}{360}$, where $\beta=(p-3)\cdot(p-5)\cdot 15-120$. Finally we take care of the solutions of the system relative to $(135)(246)$. In Lemma~\ref{counting} we have seen that it has $2\cdot(p-4)$ solutions. But we have to consider an overlapping of solutions.  Indeed, each of the $20$ subgroups  in the conjugacy class of $\langle(135)(246)\rangle$ is contained in subgroups belonging to the conjugacy class of  $\langle(123456)\rangle$. So, we have $6$ solution to take out from $2\cdot(p-4)$ and this happens for each subgroup conjugated to $\langle(135)(246)\rangle$. This leads to a total of $120$ solutions to take  out from $2\cdot(p-4)\cdot 20$. Then dividing by $240$ one gets that the number conjugacy classes whose normalizer contains a cyclic group of order $3$, which is $\frac{\delta}{240}$ where $\delta=2\cdot(p-4)\cdot 20-120$. Then an easy computation gives the rest of  the formula for case $(c)$.	
	\item[$(d)$] The proof for this case follows straightforward from those of the previous cases.  
\end{description}
\end{proof}

\begin{example}
For each of the cases described in Theorem~\ref{k5} we give an example of explicit computation of the conjugacy classes.
		\begin{table}[H]
		\centering
				\caption{Examples}
		\begin{tabular}{|c|c|c|}
			\hline 
			Prime  &   Conjugacy Classes & $m_H$ \\
			\hline
			\multirow{2}{*}{$p=17$}  &  	 $7$ & $360$	 \\  & $42$ & $720$\\
			\hline
			\multirow{3}{*}{$p=11$} &  $1$ & $144$ \\ &  $2$ & $360$ \\ &  $3$ & $720$ \\
			\hline  
			\multirow{4}{*}{$p=13$} & 	 $1$ & $120$	 \\  & 	 $3$ & $360$\\  & 	 $1$ & $240$\\
			 & 	 $9$ & $720$\\
			\hline 
			\multirow{5}{*}{$p=31$} & 	 $1$ & $120$	 \\  & 	 $1$ & $144$\\  & 	 $30$ & $360$\\
			& 	 $4$ & $240$\\ 	& 	 $775$ & $720$\\
			\hline 
		\end{tabular} 
	\end{table}
\end{example}
\section{Other admissible groups}\label{othergroups}
In this section we show that from the results for the cyclic cyclic groups previously displayed in this paper we can obtain some info regarding non-cyclic subgroups of $\mathrm{PGL}_{k+1}(\mathbb{C})$. We first recall the main classification of the finite subgroups of $\mathrm{PGL}_{k+1}(\mathbb{C})$ described by Dolgachev and Iskovskikh in~\cite{Dolgar}). This classification is borrowed from that of finite subgroups of $\mathrm{GL}_{k+1}(\mathbb{C})$, keeping in mind that each finite subgroup of $\mathrm{PGL}_{k+1}(\mathbb{C})$ can be represented by a finite subgroup of $\mathrm{GL}_{k+1}(\mathbb{C})$. Let $G$ be a finite subgroup of the general linear $\mathrm{GL}_{k+1}(\mathbb{C})$. The group $G$ is called intransitive if the representation of $G$ in $\mathbb{C}^{k+1}$ contains an invariant non-zero subspace. Otherwise it is called transitive. A transitive group $G$ is called imprimitive if it contains an intransitive normal subgroup $G_0$ . In this case $\mathbb{C}^{k+1}$ decomposes into a direct sum of $G_0$-invariant proper subspaces, and elements from $G$ permute them. A group is primitive if it is neither intransitive, nor imprimitive. The admissible
groups we are describing here are intransitive and imprimitive and are semi-direct
groups where the normal subgroup is one of the admissible groups of order $p$ we
have described in the previous section. This construction is possible  because, as we have seen in the previous section, for any $p$ there exists (up to conjugation) some group of order $p$ whose normalizer contains the group generated by an automorphism represented by a permutation matrix. We will carry out our description subdividing the discussion according to the different values of $k$ for $k\in\left\lbrace {2, 3, 4, 5}\right\rbrace $. As for notation, just for this section, given  a projective transformation
$$(x_0,\ldots, x_k)\mapsto (L_0(x_0, \ldots, x_{k}), \ldots, L_k(x_0, \ldots, x_{k}))$$
we denote it with 
$$[L_0 (x_0, \ldots, x_{k}),\ldots, L_k(x_0, \ldots, x_{k})]$$ 
and the group it generates with
$$<[L_0 (x_0, \ldots, x_{k}),\ldots, L_k(x_0, \ldots, x_{k})]>.$$
So, for example, the permutation in  $\mathrm{PGL}_{4}(\mathbb{C})$  corresponding to the the permutation $(12)(34)$ in $S_4$ will be denoted with $[x_1, x_0, x_3, x_2]$ and the admissible cyclic group of order $11$ in $\mathrm{PGL}_{3}(\mathbb{C})$ with id. vector $[2, 1, 0]^T$ will be denoted with $<[\epsilon_{11}^2x_0, \epsilon_{11}x_1, x_2]>$. Finally, the semidirect product between $2$ finite subgroups $A$ and $B$ of $\mathrm{PGL}_{k+1}(\mathbb{C})$ will be denoted with $A : B$.
\subsection{$k=2$}
The description of the groups of our interest has already been made by the  authors in~\cite{Dolgar}. Here we give some additional info and underline some important aspects.
In Theorem~\ref{1theorem} we have fully described the conjugacy classes of cyclic groups of order $p$  whose normalizer  contains the subgroup generated by an automorphism represented by a permutation matrix; elaborating this info we have then the following results. \\
Intransitive:
\begin{center}
$Int_{12}:=<[ \epsilon_{p}^{p-1} x_0 ,\epsilon_{p} x_1, x_2]> : <[x_1, x_0, x_2]>$
\end{center}
\vspace{0.3 cm}
Imprimitive:
\begin{center}
	 $Imp_{123}:=<[ \epsilon_{p}^{s} x_0 ,\epsilon_{p} x_1, x_2]> : <[x_1, x_2 , x_0]>$
\end{center}
where $s$ is a solution of $s^2-s + 1\equiv  0 \ (mod \ p)$.
\begin{center}
	 $Imp_{S_3}:=<[ \epsilon_{3}^{2} x_0 ,\epsilon_{p} x_1, x_2]> : (<[x_1, x_2 , x_0]>: <[x_1, x_0 , x_2]>).$
\end{center}

\begin{theorem}
The number of conjugacy classes of the groups $Int_{12}$,  $Imp_{123}$ and  $Imp_{S_3}$ is as follows depending on the value of $p$: 
\begin{table}[H]
	\centering
	\label{atable1}
	\caption{Number of conjugacy classes in $\mathrm{PGL}_{3}\left(\mathbb{C}\right)$}
	\renewcommand{\arraystretch}{2}
	\begin{tabular}{|c|c|c|}
		
		\hline 
		Group &	Prime Parity &  Number of Conjugacy Classes \\ 
		\hline 
		$Int_{(12)}$ &  any  & $1$ \\
		\hline 
		$Imp_{(123)}$  & $p\equiv 1 \ (mod \ 3)$ &  
		$1$ \\ 
		\hline 
		$Imp_{S_3}$  & $p=3$ &  
		$1$ \\ 
		\hline
	\end{tabular} 
\end{table}
\end{theorem}
\begin{proof}
The proof follows from the proof of Theorem 4.7 in~\cite{Dolgar} and from the results in Theorem~\ref{1theorem}.
\end{proof}

\subsection{$k=3$}

Here the data of our interest are displayed in Theorem~\ref{acase3}. For $p\geq 5$ we have then the following groups: \\
Intransitive:
\begin{center}
	$Int_{(13)(24)}:=<[\epsilon_p^s x_0, \epsilon_p^{s+1} x_1,  \epsilon_p x_2, x_3]> : <[x_2, x_3 , x_0, x_1]>$	
\end{center}
where $2\leq s\leq p-2$ is an integer.  
\begin{center}
	$Int_{(123)}:=<[\epsilon_p^{t^2} x_0, \epsilon_p^{t} x_1,  \epsilon_p x_2, x_3]> : <[x_1, x_2, x_0, x_3]>$,
\end{center}
where $t$ is a solution of $t^2-t+1\equiv 0 \ (mod \ p)$ and $p\equiv 1 \ (mod \ 3)$.\\
Inprimitive:
\begin{center}
		$Imp_{(1234)}:=<[\epsilon_p^s x_0, \epsilon_p^{s+1} x_1,  \epsilon_p x_2, x_3] : [x_1, x_2, x_3, x_0]>$
\end{center}
where $s$ is a solution of $s^2\equiv -1 \ (mod \ p)$ and $p\equiv 1 \ (mod \ 4)$.\\
For the above groups we have the following result:
\begin{theorem}
The number of conjugacy classes of the groups $Int_{(13)(24)}$,  $Int_{(123)}$ and $Imp_{(1234)}$ is as follows depending on the value of $p$:
\begin{table}[H]
	\centering
	\label{atable2}
	\caption{Number of conjugacy classes in $\mathrm{PGL}_{4}\left(\mathbb{C}\right)$}
	\renewcommand{\arraystretch}{2}
	\begin{tabular}{|c|c|c|}
		
		\hline 
		Group &	Prime Parity &  Number of Conjugacy Classes \\ 
		\hline 
		$Int_{(13)(24)}$ &  any  & $\frac{3\cdot(p-3)}{12}$ \\
		\hline 
		$Int_{(123)}$  & $p\equiv 1 \ (mod \ 3)$ &  
		$1$ \\ 
		\hline 
		$Imp_{(1234)}$  & $p=1 \ (mod \ 4)$ &  
		$1$ \\ 
		\hline
	\end{tabular} 
\end{table}
\end{theorem}
\begin{proof}
The proof follows straight from the results in Theorem~\ref{acase3} and it is similar to that of Theorem 4.7 in~\cite{Dolgar}, so we omit it.
\end{proof}
\subsection{$k=4$}
For $k=4$ the data of our interest are displayed on Table~\ref{table9} in Theorem~\ref{case4}. From these data we get the following groups:\\ 
Intransitive:
\begin{center}
	$Int_{(14)(23)}=<[\epsilon_p^{-1} x_0, \epsilon_p^{t} x_1, \epsilon_p^{-t} x_2, \epsilon_p x_3, x_4]> : <[x_3, x_2 , x_1, x_0, x_4]>$
\end{center}
where $2\leq t\leq p-1$ is an integer such that the entries of $[\epsilon_p^{-1} x_0, \epsilon_p^{t} x_1, \epsilon_p^{-t} x_2, \epsilon_p x_3, x_4]$ are all distinct,
\begin{center}
	$Int_{(1234)}:=<[\epsilon_p^{t^2} x_0, \epsilon_p^{t} x_1,  \epsilon_p^{t^3} x_2,  \epsilon_p x_3, x_4]> : <[x_1, x_3, x_0, x_2, x_4]>$,
\end{center}
where $t$ is a solution of $t^2\equiv -1 \ (mod \ p)$ and $p=1 \ (mod \ 4)$.\\
Inprimitive:
\begin{center}
	$Imp_{(12345)}:= <[\epsilon_p^s x_0, \epsilon_p^{-s^2+s} x_1,  \epsilon_p^{s^3-s^2+s} x_2, \epsilon_p x_3, x_4]> : <[x_1, x_2 , x_3, x_4,x_0]>$
\end{center}
where $s$ is solution of the $s^4-s^3+s^2-s+1\equiv 0 \ (mod \ p)$ and $p=1 \ (mod \ 5)$.\\
We list the number of conjugacy classes of the above groups in the following theorem.
\begin{theorem}
The number of conjugacy classes of the groups $Int_{(14)(23)}$,  $Int_{(1234)}$ and  $Imp_{(12345)}$ is as follows depending on the value of $p$:

\begin{table}[H]
	\centering
	\label{atable3}
	\caption{Number of conjugacy classes in $\mathrm{PGL}_{5}\left(\mathbb{C}\right)$}
	\renewcommand{\arraystretch}{2}
	\begin{tabular}{|c|c|c|}
		
		\hline 
		Group &	Prime Parity &  Number of Conjugacy Classes \\ 
		\hline 
		$Int_{(13)(24)}$ &  any  & $\frac{15\cdot(p-3)}{60}$ \\
		\hline 
		$Int_{(1234)}$  & $p\equiv 1 \ (mod \ 4)$ &  
		$1$ \\ 
		\hline 
		$Imp_{(12345)}$  & $p=1 \ (mod \ 5)$ &  
		$1$ \\ 
		\hline
	\end{tabular} 
\end{table}
\end{theorem}
\begin{proof}
The proof follows straight from the results in Theorem~\ref{case4} and it is similar to that of Theorem 4.7 in~\cite{Dolgar}, so we omit it.	
\end{proof}

\subsection{$k=5$}

For $k=5$ the data of our interest are displayed on Table~\ref{table12} and in Theorem~\ref{k5}. From these data we get the following groups:\\
Intransitive:
\small
\begin{center}
	$Int_{(14)(25)(36)}:=<[\epsilon_p^{u-v} x_0, \epsilon_p^{u-1} x_1, \epsilon_p^{u} x_2, \epsilon_p^v x_3, \epsilon_p x_4, x_5]> : <[x_3, x_4 , x_5, x_0, x_1,x_2]>$
\end{center}
\normalsize
where $2\leq u,v\leq p-1$ are integers such that the entries of $[\epsilon_p^{u-v} x_0, \epsilon_p^{u-1} x_1, \epsilon_p^{u} x_2, \epsilon_p^v x_3, \epsilon_p x_4, x_5]$ are all distinct,
\footnotesize
\begin{center}
	$Int_{(135)(246)}:=<[\epsilon_p^{i^2-i^2v} x_0, \epsilon_p^{-i^2v} x_1,  \epsilon_p^{i+v} x_2, \epsilon_p^v x_3, \epsilon_p x_4, x_5]> : <[x_2, x_3 , x_4, x_5, x_0,x_1]>$, 
\end{center}
\normalsize	
where $i$ is a primitive $3$rd root of unity modulo p and $p\equiv 1 \ (mod \ 3)$ and $2\leq v\leq p-1$ is an  integers such that the entries of $[\epsilon_p^{i^2-i^2v} x_0, \epsilon_p^{-i^2v} x_1,  \epsilon_p^{i+v} x_2, \epsilon_p^v x_3, \epsilon_p x_4, x_5]$ are all distinct,
\begin{center}
	$Int_{(12345)}:=<[\epsilon_p^{v^4} x_0, \epsilon_p^{v^3} x_1,  \epsilon_p^{v^2} x_2, \epsilon_p^{v} x_3, \epsilon_p x_4, x_5]> : <[x_1, x_2 , x_3, x_4, x_0,x_5]>$,
\end{center}
where $v$ is a solution of $v^4-v^3+v^2-v+1\equiv 0 \ (mod \ p)$ and  $p\equiv 1 \ (mod \ 5)$.\\
Inprimitive:
\scriptsize 
\begin{center}
	$Imp_{(123456)}:=[\epsilon_p^s x_0, \epsilon_p^{-s^2+s} x_1,  \epsilon_p^{s^3-s^2+s} x_2, \epsilon_p^{-s^4-s^3+s^2-s} x_3, \epsilon_p x_4, x_5] : [x_1, x_2 , x_3, x_4,x_5,x_0]$
\end{center}
\normalsize
where $s$ is solution of $s^2+s+1\equiv 0 \ (mod \ p)$ 
and $p\equiv 1 \ (mod \ 3)$.\\
We list the number of conjugacy classes of the above groups in the following theorem.
\begin{theorem}
The number of conjugacy classes of the groups $Int_{(14)(25)(36)}$,  $Int_{(135)(246)}$, $Int_{(12345)}$ and  $Imp_{(123456)}$ is as follows depending on the value of $p$:
\begin{table}[H]
	\centering
	\label{atable3}
	\caption{Number of conjugacy classes in $\mathrm{PGL}_{5}\left(\mathbb{C}\right)$}
	\renewcommand{\arraystretch}{2}
	\begin{tabular}{|c|c|c|}
		
		\hline 
		Group &	Prime Parity &  Number of Conjugacy Classes \\ 
		\hline 
		$Int_{(14)(25)(36)}$ &  any  & $\frac{15\cdot(p-3)\cdot(p-5)}{360}$ \\
		\hline 
		$Int_{(135)(246)}$  & $p\equiv 1 \ (mod \ 3)$ &  
		$\frac{2\cdot (p-4)\cdot 20}{240}$ \\ 
		\hline 
		$Int_{(12345)}$  & $p=1 \ (mod \ 5)$ &  
		$1$ \\ 
		\hline
			$Imp_{(123456)}$  & $p=1 \ (mod \ 3)$ &  
		$1$ \\
			\hline 
	\end{tabular} 
\end{table}
\end{theorem}
\begin{proof}
The proof follows straight from the results in Theorem~\ref{k5} and it is similar to that of Theorem 4.7 in~\cite{Dolgar}, so we omit it.
\end{proof}

\section{Cyclic groups of order $p$ with $p< k+1$}\label{power}
In this section we consider the case  $p<k+1$.  In~\cite{Mari2} we have described the necessary and sufficient  conditions for a diagonal matrix of order $p<k+1$ to represent a generator  of an admissible subgroup  $C_p$ of order $p$ in $\mathrm{PGL}_{k+1}(\mathbb{C})$. More specifically, we have shown  that the configuration of the diagonal entries depends on the $p$-parity of $k+1$.   Indeed, the $p$-parity of $k+1$ determines an integer $l$ as follows
\begin{description}
	\item if $k+1\equiv 0 \ (mod \ p)$ then $l=\frac{k+1}{p}$.
	\item if $k+1\equiv a \ (mod \ p)$ where $1 \leq a \leq p-1$, then $l=\frac{k+1-a}{p}+1$.
\end{description}
If  $k+1\equiv 0 \ (mod \ p)$ then the integer $l$ gives the size of each block of eigenvalues on the diagonal of the matrix of order $p$ which represents a generator of the group: there are $p$-blocks of size $l$ one for each power $\epsilon_p^{\alpha}$ (for $\alpha =0,1,\ldots, p-1$)  of $\epsilon_p$ where $\epsilon_p$ is a primitive $p$-th root of unity. While if $k+1\equiv a \ (mod \ p)$ on the diagonal we have $a$ blocks of size $l=\frac{k+1-a}{p}+1$ and $p-a$ blocks of size $l-1$. 
\begin{example}
	If one  take $k=5$ and  $p=3$, since $k+1\equiv 0 \ (mod \ 3)$ there are $3$ blocks of the  same size which is $2$. Following the convention introduced in 	Remark~\ref{terminology}  a generator of this group can be represented by the following vector

	\begin{equation}\label{eqn:blockmatrix10}
	\left[
	\begin{array}{c}
	2 \\
	2 \\
	1\\ 
	1\\
	0 \\ 
	0\\
	\end{array}\right]
	\end{equation}
	On the other hand for $k=4$ and $p=3$ we have $k+1\equiv 2 \ (mod \ 2)$  and so $3$ blocks:  $2$ of size $2$ and one of size $1$.
	Here again a generator can be represented by the following vector

	\begin{equation}\label{eqn:blockmatrix11}
	\left[
	\begin{array}{c}
	2 \\
	2 \\
	1\\ 
	1\\
	0 \\ 
	\end{array}\right]
	\end{equation}
\end{example}

In~\cite{Mari2} we have also  computed the number of admissible non-trivial projective representations of an abstract group $\hat{C}_p$ of order $p$, where for admissible we intend those  representations whose  image can be the automorphism group of  a  point set. We report this result for the case $p<k+1$ as we will use it here for our purposes. For more reference on the projective representations see~\cite{GKar}.  

\begin{theorem}\label{representation}
	Let  $p<k+1$ be  prime. The number $e$ of equivalence classes of non-trivial admissible  projective representations of $\hat{C}_p$ of degree $k+1$ is as follows:
	\begin{itemize}
		\item  if $k+1\equiv 0 \ (mod \ p)$  then $e=1$;
		\item  if $k+1\equiv a \ (mod \ p)$ with $1\leq a \leq p-1$ then $e=\frac{1}{a}\cdot\binom{p-1}{a-1}$.
	\end{itemize}
\end{theorem} 
Now Theorem~\ref{representation}  helps us to count straightforward  the number of conjugacy classes of admissible cycle  subgroups of order $p$ $\mathrm{PGL}_{k+1}(\mathbb{C})$ for some specific cases. Indeed, it  is clear that when the number $e$ of admissible representation is equal to  $1$, there is also one conjugacy class of subgroups. Note that this covers in particular the cases $a=1$ and $a=p-1$. We formalize these remarks in the following Corollary
\begin{corollary}\label{representation1}
	Let  $p<k+1$ be  prime.  Suppose that one of the following  conditions is fulfilled
	\begin{description}
		\item [$(a)$] $k+1\equiv 0 \ (mod \ p)$;
		\item [$(b)$] $k+1\equiv a \ (mod \ p)$ with $a\in\lbrace   {1,p-1}\rbrace$.
	\end{description}
	Then, there is only one conjugacy class of admissible groups of order $p$ in $\mathrm{PGL}_{k+1}(\mathbb{C})$.	
\end{corollary}
\begin{remark}
	Note that Corollary~\ref{representation1} tells us in particular that when $p<k+1$  and $p=2$ or $p=3$ there is always only one conjugacy class of admissible subgroups of order $p$.
\end{remark}

In order to treat the cases left out by Corollary~\ref{representation1}  we need to further exploit the idea in the argument of the proof for case $(b)$ in the same theorem. This idea is that it is possible to restrict the counting to one  of the two parts  in which the entries on the diagonal matrix are split by the integer $a$.  This generalization is done in the following theorem

\begin{theorem}\label{aconjugacy}
	Let $p\geq 5$ be a prime with $p<k+1$. Suppose further that $k+1\equiv a \ (mod \ p)$ with $2\leq a \leq p-2$.  Then, the number of conjugacy classes of  admissible subgroups of order $p$ in $\mathrm{PGL}_{k+1}(\mathbb{C})$ is the same as the number of conjugacy classes of  admissible subgroups of order $p$  in $\mathrm{PGL}_{a}(\mathbb{C})$ and in $\mathrm{PGL}_{p-a}(\mathbb{C})$.
\end{theorem}
\begin{proof}
	Let $p$ be a prime with $p<k+1$. Suppose further that $k+1\equiv a \ (mod \ p)$ with $2\leq a \leq p-2$.  
	Up to  rescaling  we can assume  that the diagonal matrix of finite order $p$ representing  a generator of an admissible subgroup of order $p$ in $\mathrm{PGL}_{k+1}(\mathbb{C})$ contains the eigenvalue $1$ at the last  of the first $a$ blocks of eigenvalues.  Moreover, up to replacing this matrix with a suitable power of  itself, we can suppose that second to the last block in the first $a$ one contains the eigenvalue $\epsilon_p$. It is clear that the choice of the eigenvalues for the remaining first $a-2$ blocks determines the set of eigenvalues in the last $p-a$ blocks. Moreover, the  last $p-a$ blocks are given up to conjugation and so  we can disregard the order of the last $p-a$ eigenvalues. Now, considering each block of  eigenvalues in $a$  as a single  eigenvalue allows to reduce the counting to the counting of conjugacy classes of  cyclic groups of order $p$ in $\mathrm{PGL}_{a}(\mathbb{C})$. Then exchanging the roles of $a$ and $p-a$ ones proves the  second part of  the statement in Lemma~\ref{aconjugacy}. 
\end{proof}
Theorem~\ref{aconjugacy} has three applications in this paper: two will be discussed here below and the third one in the next subsection .  The first application  is the computation of the number of admissible conjugacy classes of subgroups of order $p$ for $p<k+1$  using the knowledge of the number of conjugacy classes of order $p$ in $\mathrm{PGL}_{a}(\mathbb{C})$ or in $\mathrm{PGL}_{p-a}(\mathbb{C})$.  The second application is  to  compute the number  of admissible conjugacy classes in the projective linear groups of higher degree $\mathrm{PGL}_{p-a}(\mathbb{C})$ knowing that in $\mathrm{PGL}_{a}(\mathbb{C})$. Both applications rely on the results  we  gave in Section~\ref{section1}  for $p\geq k+1$ and  $3\leq k+1\leq 6$.  Thus the range of the cases covered has a limitation  due  to the limited range of the original results. Anyway, Theorem~\ref{aconjugacy}  allows us to fully compute the number of conjugacy classes of subgroups order $p$ with $p<k+1$ whenever  at least one of the $2$ integers $a$ and $p-a$  stays in the set $\lbrace {2,3,4,5,6}\rbrace$. The case where at least one between $a$ and $p-a$ is equal to $2$ is not comprised in the results of Section~\ref{section1}. However, it is a well known fact that  there is just one conjugacy class in $\mathrm{PGL}_{2}(\mathbb{C})$ of subgroups of order $p$. So also these values of $a$ are comprised.  In addition, note that the cases covered by Theorem~\ref{aconjugacy} together with those of Corollary~\ref{representation1} permit the full treatment of the problem when $p<k+1$ and $p$ lies in the set $\lbrace {5,7,11,13}\rbrace$. The data regarding these primes  are listed in the following  table

\begin{table}[H]
	\label{table15}
	\centering
	\caption{}
	\renewcommand{\arraystretch}{1.3}
	\begin{tabular}{|c|c|c|}
		\hline 
		Prime & Values of $a$	&  Number of Conjugacy Classes \\
		\hline
		$p=5$ &  $0, 1, 2, 3, 4$ & $1$ \\ 
		\hline
		
		\multirow{2}{*}{$p=7$} & $0,1,2,5,6$ & 	 $1$ \\ 	
		& $3,4$ & $2$ \\
		\hline
		\multirow{4}{*}{$p=11$} & $0,1,2,9,10$ & 	 $1$ \\ 	
		& $3,8$ & $2$ \\
		& $4,7$ & $4$ \\
		& $5,6$ & $6$ \\
		\hline
		\multirow{5}{*}{$p=13$}  &  $0,1,2,11,12$		 &  $1$	 \\
		& $3,10$ &  $3$ \\ 
		& $4,9$ &  $7$ \\ 
		& $5,8$ &  $10$ \\ 
		& $6,7$ &  $14$ \\ 	
		\hline
	\end{tabular} 
\end{table}  
Finally, we consider the second application of Theorem~\ref{aconjugacy}  that is the direct exploitation of the kind of  duality between $\mathrm{PGL}_{a}(\mathbb{C})$ or in $\mathrm{PGL}_{p-a}(\mathbb{C})$. This duality allows us to treat the case $p\geq k+1$  in projective linear groups of degree higher than $6$ for some specific pair of values of $p$ and $k+1$. Indeed by means of the results in Section $3$, those in Theorem~\ref{aconjugacy} and in Corollary~\ref{representation1}  we can compute for each prime $p$ the number the conjugacy classes in $\mathrm{PGL}_{p-a}(\mathbb{C})$ when $a$ ranges in the set  $\lbrace {1,2,3,4,5,6}\rbrace$.  Note that it is possible to  add to these cases also the case $k+1=p$ as it is clear that here there is only one conjugacy class. As an example  of what we have just discussed we list in the table below the number of conjugacy classes for the prime $p=19$ for $k+1$ assuming values in the set $\lbrace {19,18,17,16,15,14, 13}\rbrace$.

\begin{table}[H]
	\label{table16}
	\centering
	\caption{Case $p=19$}
	\renewcommand{\arraystretch}{1.3}
	\begin{tabular}{|c|c|}
		\hline 
		$k+1$ 	&  Number of Conjugacy Classes \\
		\hline
		$19$ & 1 \\
		\hline 
		$18$ & 1 \\
		\hline
		$17$ & 1 \\
		\hline
		$16$ & 4 \\
		\hline
		$15$ & 14 \\
		\hline
		$14$ & 36 \\
		\hline
		$13$ & 86 \\
		\hline
	\end{tabular} 
\end{table}

\begin{remark}\label{groupassociation}
	In the rest of this paper  we will refer to the relation described in Theorem~\ref{aconjugacy} between  a subgroup $H_1$  of $\mathrm{PGL}_{a}(\mathbb{C})$  of order $p$ and a subgroup $H_2$  of $\mathrm{PGL}_{p-a}(\mathbb{C})$ of order $p$  as association and the two subgroups will  be referred to as associated subgroups. 
\end{remark}

 \subsection{Relation with the association between point sets}\label{galetransform}
In this subsection we intend to show  a beautiful relation connecting the association between cyclic subgroups of order $p$ described in the proof of Theorem~\ref{aconjugacy} (see also Remark~\ref{groupassociation}) and the association between point sets.
We start with a description of the association between sets of points. For more reference on this concept see~\cite{Cobl}, ~\cite{Cobl1}, ~\cite{Dolga} and ~\cite{Eise}.\\

In this subsection $\mathbb{P}^k(\mathbb{C})^{n}$ is $n$-th  product  of $\mathbb{P}^k(\mathbb{C})$ with itself. 
\begin{definition}
	Let $A\in (\mathbb{P}^k)^n(\mathbb{C})$ be a set consisting of $n$ ordered points  $P_i$. Choosing projective coordinates  $[\alpha_{i 0}, \ldots, \alpha_{ik}]$ for each $P_i$, the matrix $M_A$ of the projective coordinates of $A$ is the $(k+1)\times n$ matrix whose columns are the projective coordinates of the points $P_i$, i.e.
	
	\begin{center}
		$M_A=\left(\begin{array}{lll}
		\alpha_{1 0} & \cdots & \alpha_{n 0}\\
		\vdots &  & \vdots\\
		\alpha_{1 k} & \cdots  & \alpha_{n k}
		\end{array}\right)$.
	\end{center}
\end{definition}

\begin{remark}
	Due to the choice of the projective coordinates of the points $P_i$, the matrix $M_A$ is given up to multiplication of each column by a constant.
\end{remark}

\begin{example}\label{association2}
	Let $A\in (\mathbb{P}^1)^5(\mathbb{C})$ be the ordered set consisting of the points\\
	
	$P_1=\left[1,0\right]$, $P_2=\left[0,1\right]$, $P_3=\left[1,1\right]$, $P_4=\left[1,\sqrt{2}\right]$, $P_5=\left[1,-\sqrt{2}\right]$. \\
	
	Then the matrix of the projective coordinates of the points in $A$ is \\
	
	\begin{center}
		$M_A=\left(
		\begin{array}{ccccc}
		1 & 0 & 1 & 1 & 1\\
		0 & 1 & 1 & \sqrt{2} & -\sqrt{2} 
		\end{array}\right)$.
	\end{center}
\end{example}

Now we restrict our attention to general sets of points. Recall that a set of points $A\in (\mathbb{P}^k)^n(\mathbb{C})$ is said to be general if any subset of $t\leq k+1$ points spans a $(t-1)$-dimensional linear projective space. Hereafter in this section, we denote the subset of the general point sets in $(\mathbb{P}^k)^n(\mathbb{C})$ with $\Gamma^k_n(\mathbb{C})$.  

\begin{definition}\label{associatepoint}
	Let $A\in \Gamma^k_n(\mathbb{C})$ and let $M_A$ be its matrix of projective coordinates, a set $B\in \Gamma^{n-k-2}(\mathbb{C})$ is said to be associated to $A$ if its coordinate matrix $M_B$ satisfies
	
	\begin{center}
		$M_A\cdot\varLambda\cdot {}^t M_B =0$ \end{center}
	for some diagonal matrix $\varLambda=diag(\delta_1, \ldots, \delta_n)$ with all $\delta_i\neq 0$.  
\end{definition}

\begin{example}
	The general set of points in Example~\ref{association2}  is associated to the set $V\in\Gamma^2_5(\mathbb{C})$ consisting of the points
	\begin{center}
		
		$Q_1=\left[1, 1, 1\right]$, $Q_2=\left[1, \sqrt{2}, -\sqrt{2}\right]$, $Q_3=\left[1, 0, 0\right]$, $Q_4=\left[0, 1, 0\right]$, $Q_5=\left[0, 0, 1\right]$.                                                                                                                                                          \end{center}
	
	Indeed, taking the $5\times 5$ matrix $\varLambda=diag(1,1,-1,-1,-1)$ one sees that\begin{center}

		$M_A\cdot\varLambda\cdot {}^t M_B =0$                                     \end{center}
\end{example}

\begin{remark}\label{permute}
	The relation of association is symmetric. Note also that for every permutation $\pi\in S_n$ the sets $A=\{P_1,\ldots, P_n\}$ and $B=\{Q_1,\ldots, Q_n\}$ are
	associated if and only if $\lbrace P_{\pi(1)},\ldots, P_{\pi(n)}\rbrace$ and $V=\lbrace Q_{\pi(1)},\ldots, Q_{\pi(n)}\rbrace$ are associated.
\end{remark}

The association between point sets determines an isomorphism between the orbits of the associated points, i.e.

\begin{center}
	$\varphi_{k,n-k-2}: \Gamma^k_n(\mathbb{C})/\mathrm{PGL}_{k+1}(\mathbb{C})\rightarrow\Gamma^{n-k-2}_n(\mathbb{C})/\mathrm{PGL}_{n-k-1}(\mathbb{C})$.                                                                                                                  \end{center}

Unless $k=n-k-2$, the associated sets of points are in spaces of different dimension. As pointed out  by Coble in~\cite{Cobl}, the conventional methods of passing from one space to another (and from one point set to a point set associated to it) are the process of mapping the space of lower dimension to that of higher dimension and the process of projecting from the space of higher dimension upon that of lower dimension. Let us consider for example the mapping from $\mathbb{P}^1(\mathbb{C})$. The linear system $O_{\mathbb{P}^1}(n-3)$ of hypersurfaces of degree $n-3$ defines an embedding of $\mathbb{P}^1$ in $\mathbb{P}^{n-3}$ of degree $n-3$ (the Veronese embedding), whose image is a rational normal curve. This embedding maps a set of points in $\Gamma^1_n(\mathbb{C})$ to its associate in $\Gamma^{n-3}_n(\mathbb{C})$. 
There are actually other methods for finding a point set associated to  a given point set. One of these is described in~\cite{Eise}. Given a point set $T$ in $\mathbb{P}^k(\mathbb{C})$ of  cardinality $n$, one replaces it, if necessary, with a point set in the same projective class whose first $k+1$ points are the fundamental points in $\mathbb{P}^k(\mathbb{C})$  taken in the canonical  order. A matrix of the coordinates for  $T$ has then the form
\begin{center}
	$(I_{k+1}|A)$
\end{center}  
Then it is easily to verify that the point set represented by the matrix
\begin{center}
	$(A^T|I_{n-k-2})$
\end{center}
is associated to $T$. In the last part of this section we will show a easier method in the special case of the point sets with cyclic automorphism of order $p$.\\
The definition of association between ordered point sets easily extends to unordered point sets.

\begin{definition}\label{association23}
	A point set $T\in\mathbb{P}^k_n(\mathbb{C})$ and a point set $V\in\mathbb{P}^{n-k-2}_n(\mathbb{C})$  are said to be associated if there exists an ordering of the points in $T$ and an ordering of the points in $V$ such the corresponding ordered sets are associated in the sense of Definition~\ref{associatepoint}.
\end{definition}

It is easy to prove that associated point sets have isomorphic automorphism groups

\begin{lemma}\label{associatedgroups}
	Two associated point-sets have isomorphic automorphism groups.
\end{lemma}
\begin{proof}
	Let  $T$ be a point set in $\mathbb{P}^k_n(\mathbb{C})$ and let $V\in \mathbb{P}^{n-k-2}_n(\mathbb{C})$ be any point set associated to $T$.  By definition of association between unordered sets of points, there exists an ordering of the points in $T$ and in $V$ such the corresponding ordered sets are associated. Let $T^o$ and $V^o$ the ordered point sets corresponding to these orderings, then by the definition of association
	\begin{center}
		$M_{T^o}\cdot\varLambda\cdot {}^t M_{V^o} =0$, \end{center}
	where $\varLambda$ is a $n\times n$ invertible diagonal matrix.
	Let $f$ be an automorphism of $T$. The automorphism $f$ permutes the points in $T$, so we have 
	
	\begin{center}
		$M_{(T^o)^f}= M_{T^o}\cdot P_f$
	\end{center}
	where $P_f$ is an $n\times n$ permutation matrix.
	Since $M_{(T^o)^f}$ also satisfies
	
	\begin{center}
		$M_{(T^o)^f}\cdot\varLambda\cdot {}^t M_{V^o} =0$. \end{center}
	we get
	\begin{equation}\label{identity1}
	(M_{T^o}\cdot P_f) \cdot\varLambda\cdot {}^t M_{V^o} =0. \end{equation}
	On the other hand, as it has been explained in Remark~\ref{permute},
	\begin{equation}\label{identity2}
	(M_{T^o}\cdot P_f) \cdot\varLambda\cdot {}^t (M_{V^o}\cdot P_f) =0. \end{equation}
	Let $(V^o)'$ be the ordered set of points corresponding to the matrix $M_{V^o}\cdot P_f$. Identities~\ref{identity1} and~\ref{identity2} show that $V^o$ and $(V^o)'$ stay in the same $\mathrm{PGL}_{n-k-1}\left(\mathbb{C}\right)$-orbit, that is there exists $g\in \mathrm{PGL}_{n-k-1}\left(\mathbb{C}\right)$ such that $(V^o)'=(V^o)^g$. Note that $g$ is uniquely determined by the permutation $P_f$ as $n\geq k+3$.  Thus we have shown the existence of a map $\varphi: \mathrm{Aut}(T)\rightarrow \mathrm{Aut}(V)$. This map is clearly injective and an homomorphism. Exchanging the roles of $T$ and $V$ one then sees that this map is also onto.      \\
\end{proof}

\begin{theorem}\label{correspondance}
	Let $p\geq 5$ be a prime. Suppose further that $a$ is an integer with  $2\leq a \leq p-2$.  Let $H$ be an  admissible subgroup of order $p$ in $\mathrm{PGL}_{a}(\mathbb{C})$ in diagonal form and let $K$ be one of its associated subgroups in $\mathrm{PGL}_{p-a}(\mathbb{C})$ described in Theorem~\ref{aconjugacy} (see Remark~\ref{groupassociation}).
	Let $T$ be the point set in $\mathbb{P}^{a-1}\left(\mathbb{C}\right)$ consisting of the orbit of the identity point under the action of the group $H$. Similarly, let $V$ be the point set in $\mathbb{P}^{p-a-1}\left(\mathbb{C}\right)$ consisting of the orbit of the identity point under the action of the group $K$. Then the point sets $T$ and $V$ are associated.
\end{theorem}
\begin{proof}
	In order to make things easier we can reason in terms of vectors for the generators of the groups $H$ and $K$. For the group $H$ we take the id. vector and we denote it with $h$.  For the second  group we chose any vector whose entries are the complement in $\mathbb{Z}_p$ of the entries of $h$ and denote it with $k$. Now, let $T$ and $V$ be as in  the statement of Theorem~\ref{correspondance}. We fix an ordering of $T$ such that the first point is  represented by the vector $h$ and  the other points are in sequence represented by the increasing multiples of $h$. We do the same  for  $V$ starting by the point represented by the vector $-k$. Let us indicate with $A$ and $B$ the matrices, whose columns are the representative so chosen for each point in the now ordered sets $T$ and $V$.  By Definitions~\ref{associatepoint} and~\ref{association23} we need to  show that the product $A\cdot B^T$ is a null $(a\times (p-a))$ matrix. First note  that  it is enough to show  that the first term in each product row by column is not null. Indeed, if this is the case then all  the terms appearing in a single  row by column product would be all distinct, as they are all the multiples of a non-null element of $\mathbb{Z}_p$. Reasoning in terms of corresponding $p$-th roots of unity
	this shows that the terms in each row by column product run all over the $p$-th roots. Hence  their sum is  zero. So, to prove Theorem~\ref{correspondance} we need to prove the condition for the first terms. We have to reason on the entries of  $k$ and $-k$. If $k$ contains the opposite of any entry of $h$, then taking $-k$ eliminates the problem. On the other hand if $k$ contains no opposite of elements in $h$, so does $-k$ as the entries of $k$ are taken from the complement in $\mathbb{Z}_p$ of those of $h$. So we are done with the proof of the statement in Theorem~\ref{correspondance}.
\end{proof}

\begin{example}\label{grouppoint}
	Let us consider the cyclic group $H$ of order $11$ in $\mathrm{PGL}_{4}(\mathbb{C})$  with id. vector $[4,2,1,0]$. Consider now the associated cyclic group $K$ in  $\mathrm{PGL}_{7}(\mathbb{C})$,  with vector $[10,9,8,7,6,5,3]$ (recall that the entries of this  vector are taken from the complement in $\mathbb{Z}_{11}$ of the entries of the id. vector for $H$). Consider then the two $11$-point sets $V$ and $T$ consisting respectively of the orbits of the identity points in   $\mathbb{P}^3(\mathbb{C})$ and $\mathbb{P}^6(\mathbb{C})$  under the action of the groups $H$ and $K$. We want to show that $T$ and $V$ are associated.
	Let us now fix an order on $T$ such that the first point is represented by the vector $(4,2,1,0)$ and the other in sequence by the multiple of this one and  an analogue  order on $V$ with the first point represented by the vector $(1,2,3,4,5,6,8)$. Denoting with $A$ and $B$ the corresponding matrices for the point sets  $T$ and $V$, one gets
	\begin{center}
		$A=\left(\begin{array}{ccccccccccc}
		4 & 8 & 1 & 5 & 9 & 2 & 6 & 10 & 3 & 7 & 0 \\ 
		2 & 4 & 6 & 8 & 10 & 1 & 3 & 5 & 7 & 9 & 0 \\ 
		1 & 2 & 3 & 4 & 5 & 6 & 7 & 8 & 9 & 10 & 0 \\ 
		0 & 0 & 0 & 0 & 0 & 0 & 0 & 0 & 0 & 0 & 0
		\end{array}\right) $
	\end{center}
	and
	\begin{center}
		$B^T=\left(\begin{array}{ccccccc}
		1 & 2 & 3 & 4 & 5 & 6 & 8 \\ 
		2 & 4 & 6 & 8 & 10 & 1 & 5 \\ 
		3 & 6 & 9 & 1 & 4 & 7 & 2 \\ 
		4 & 8 & 1 & 5 & 9 & 2 & 10 \\ 
		5 & 10 & 4 & 9 & 3 & 8 & 7 \\ 
		6 & 1 & 7 & 2 & 8 & 3 & 4 \\ 
		7 & 3 & 10 & 6 & 2 & 9 & 1 \\ 
		8 & 5 & 2 & 10 & 7 & 4 & 9 \\ 
		9 & 7 & 5 & 3 & 1 & 10 & 6 \\ 
		10 & 9 & 8 & 7 & 6 & 5 & 3 \\ 
		0 & 0 & 0 & 0 & 0 & 0 & 0
		\end{array} \right) $
	\end{center}
Now, reinterpreting the entries of $A$ and $B$ as powers of a primitive $11$-th root of unity $\epsilon_{11}$, it is easy to see that the product of these $2$  matrices is a null matrix, as the first term in each row by column product is a power of $\epsilon_{11}$ different from $1$. Hence the point sets are associated. 
\end{example}
Starting from the result given in Theorem~\ref{correspondance} it is possible to extend the method to any union of non trivial orbits of an admissible cyclic group $C_p$. Here it is exploited the interplay between admissible cyclic groups of order $p$ in $\mathrm{PGL}_{a}(\mathbb{C})$ (or in $\mathrm{PGL}_{p-a}(\mathbb{C})$) and those in $\mathrm{PGL}_{k+1}(\mathbb{C})$ when $k+1\equiv a \ (mod \ p)$.
\begin{example}\label{generalization}
	In this example in order to avoid confusion with the coordinates of the points involved, we will describe the generators of the groups with the id. matrices.
	Let us then consider again the cyclic group $H$ in Example~\ref{grouppoint}. Now we consider the point set $T_1$ in $\mathbb{P}^3(\mathbb{C})$ consisting of the union of the point set $T$ (the orbit of the point $P_1=[1,1,1,1]$) and the orbit of the point $P_2=[4,3,2,1]$ under the same group. We already fix an order on $T_1$, taking first the points in $T$ in the order described in Example~\ref{grouppoint} followed by the points in the orbit of $P_2$  taken in the analogous order. We also fix the representatives of the coordinates of these points in the usual manner  and denote the matrix of the coordinates with $A_1$. Now,  the cardinality of $T_1$ is $22$, so any point set associated  to $T_1$ stays in  $\mathbb{P}^{17}(\mathbb{C})$. Now, since $18\equiv 7 \ (mod \ 11)$  any admissible subgroup in $\mathrm{PGL}_{18}(\mathbb{C})$ has, in terms of multiplicity of the eigenvalues, the shape we have described at the beginning of this section.  This gives us a hint on how to  build a point set associated to $T_1$. We may assume that the entries of the id. matrix of the group in question consist  of those of the group $H$ in Example~\ref{grouppoint} taken twice plus the remaining $4$ given by the inverse of each diagonal entry of $H$. Moreover, we can choose the order for these $4$ entries so that the are in sequence the inverse of those of $H$.  In practice we choose as id. matrix for this group the following matrix $E$
	\begin{center}
		$E=diag(a,a^2,a^3,a^4,a^5,a^6,a^8,a,a^2,a^3,a^4,a^5,a^6,a^8,a^7,a^9,a^{10},1)$.
	\end{center}
	where $a=\epsilon_{11}$ is the same  primitive $11$-th root of unity chosen for $H$. Let us denote with $K_1$ the group generated by this matrix.
	Now, we build an associated point set taking $2$ non-trivial orbits under the action of $K_1$. These orbits are respectively that of the point \begin{center}
		$Q_1=[1,1,1,1,1,1,1,0,0,0,0,0,0,0,-4,-3,-2,-1]$
	\end{center} and that of the point
	\begin{center}
		$Q_2=[0,0,0,0,0,0,0,1,1,1,1,1,1,1,1,1,1,1]$
	\end{center}
	taken in this order.  Moreover, in each orbit we take as first point $EQ_i$ for $i\in \langle{1,2}\rangle$ followed by its conjugates according to the sequence of powers of $E$. Let us denote with $B_1$ be the matrix of coordinates corresponding to this point set. It is then easy to check that $A_1\cdot B_1^T$ is a null matrix. Indeed, the product of the rows of $A_1$ by each of the first $14$ columns of $B_1^T$   gives always $0$ for what we have seen in Example~\ref{grouppoint}. We are left with the products of the rows of $A_1$ with  the last $4$ columns of $B_1^T$. It is enough to consider the product of each of the rows of $A_1$ with the $15$-th columns of $B_1^T$ has for the other columns the argument is similar. Let us start with the first row of $A_1$.  The first element in this product is now $-4$ and so are all the others till the $7$-th row. From the $8$-th row  they are all equal to $4$. Summing everything up this gives $0$. On the other hand for each of the rows from $2$ to $4$ in $A$, the product with the $15$-th column, gives two sums of $7$ products and each of these sums is equal to $0$. What seen for the $15$-th column clearly applies to  the remaining columns. So, the two point sets are associated.
	The procedure can be iterated to three non-trivial orbits of $H$, adding to the points in $T_1$ for example the orbit of the point $P_3=[7,5,3,1]$. We get a point set $T_2$ on which we fix the order analogous to that on $T_1$. We also fix the representatives of the coordinates of these points in the usual manner  and denote the matrix of the coordinates with $A_2$.
	Now, $T_2$ has cardinality $33$, so  we are looking for a point set in  $\mathbb{P}^{28}(\mathbb{C})$. Since $29\equiv 8 \ (mod \ 11)$ and $\frac{22}{11}=2$ this leads to think that the suitable cyclic group in  $\mathrm{PGL}_{29}(\mathbb{C})$ has to be represented by the id. matrix
	\begin{center}
		$G=diag(\bar{a},\bar{a},\bar{a},\bar{b},\bar{b})$.
	\end{center}
	where $\bar{a}$  stands for the sequence of entries $a, a^2, a^3, a^4, a^5, a^6, a^8$ and  
	$\bar{b}$ for the sequence of entries $a^7, a^9, a^{10}, 1$. Let us denote with $K_2$ this group. Now to build our associated point set we need to find three suitable non-trivial orbits for $K_2$. Reasoning as for the previous case it is easy to see that it is enough to consider the orbits of the following three points
	
	\footnotesize
	\begin{center}
		$R_1=[1,1,1,1,1,1,1,0,0,0,0,0,0,0,0,0,0,0,0,0,0,-4,-3,-2,-1,-7,-5,-3,-1]$\\ \vspace{0.1 cm}	
		$R_2=[0,0,0,0,0,0,0,1,1,1,1,1,1,1,0,0,0,0,0,0,0,1,1,1,1,0,0,0,0]$ \\ \vspace{0.2 cm}
		$R_3=[0,0,0,0,0,0,0,0,0,0,0,0,0,0,1,1,1,1,1,1,1,0,0,0,0,1,1,1,1]$ \\
	\end{center}
	
	\normalsize 
	Let us call this point set  $V_2$. Now, let us fix on  $V_2$ an order analogous to that on $T_2$ with the orbits of the points $R_1$, $R_2$ and $R_3$ taken in this sequence. We also choose in an  analogous manner the representatives of the coordinates one each point in $V_2$. Then denoting with  $B_2$ the matrix of coordinates so obtained,  one easily sees that $A_2\cdot B_2^T=0$ and that the two point sets $T_2$ and $V_2$ are associated.
\end{example}
The method described in Example~\ref{generalization} can be generalized  to any admissible cyclic group of order $p$ and to any of union of its non-trivial. We summarize
it here below the procedure. 

\begin{enumerate}
	\item The starting point is an admissible group $H$ of order $p$ in $\mathrm{PGL}_{a}(\mathbb{C})$ (where $2\leq a\leq p-2$) with id vector $h$.
	\item The second step is to take $l\geq 1 $  non-trivial orbits of this group, to get a point set of order $l\cdot p$ in $\mathbb{P}^{a-1}(\mathbb{C})$, with the first orbit to be that of the identity point.
	\item The third step is to individuate a corresponding admissible subgroup in  $\mathrm{PGL}_{lp-a}(\mathbb{C})$ following the description given at the beginning of this section. This is done in practice as follows:
	\begin{description}
		\item [$(a)$] Consider the complement  $k$ of $h$ in $\mathbb{Z}_p$. Let denote with $-k$  the subset of $\mathbb{Z}_p$ consisting of the opposite of the elements in $k$.  
		\item [$(b)$] Compute $u=\frac{lp-a}{p}$ and $t=u+1$ and take $t$ copies of $-k$. 
		\item [$(c)$] Take the complement of $-k$ in $\mathbb{Z}_p$ and order the elements in this set so that they are in sequence the opposite of the entries in $v$. Let us denote this set with $z$ .
		\item [$(d)$] Consider the group $K_{l-1}$ in $\mathrm{PGL}_{lp-a}(\mathbb{C})$ represented by the vector whose entries are in sequence the elements in the $t$ copies (copy after copy) of $k$ and $u$ copies of $z$.
	\end{description}
	\item Finally take the $l$ non-trivial orbits of $K_{l-1}$ in $\mathbb{P}^{lp-a-1}(\mathbb{C})$ in the way described in Example~\ref{generalization}.
\end{enumerate} 

Clearly the method can be reversed starting from a point set  in $\mathbb{P}^{lp-a-1}(\mathbb{C})$ of the shape described  in Example~\ref{generalization} and going back to a point set in   $\mathbb{P}^{a-1}(\mathbb{C})$.

\begin{remark}
	Note that the peculiar configuration of the points in $\mathbb{P}^{lp-a-1}(\mathbb{C})$ described above is consistent with the normal form we described in~\cite{Mari2} for a point set consisting of $l$ non-trivial orbits of an admissible cyclic group of order $p$ in $\mathrm{PGL}_{lp-a}(\mathbb{C})$. Indeed the group described in the method above is a conjugate of that described in~\cite{Mari2}. 
\end{remark}

\end{document}